\documentclass[12pt]{article}

\usepackage{amsmath}
\usepackage{amsthm}
\usepackage{amssymb}
\usepackage{cite}
\usepackage{url}
\usepackage[multiple]{footmisc}
\usepackage{graphicx}
\usepackage{epstopdf}
\usepackage{chngcntr}
\usepackage{enumitem}
\usepackage{hhline}
\usepackage{array}
\usepackage{hyperref}
\usepackage{color}

\usepackage{makecell}
\usepackage{multirow}

\hypersetup{colorlinks=false}

\if@twoside
\else 
\fi

\numberwithin{equation}{section}
\counterwithin{figure}{section}
\counterwithin{table}{section}

\theoremstyle{plain}

\theoremstyle{definition}

\theoremstyle{remark}
\newtheorem*{remark}{Remark}

\newcommand{\norm}[1]{\left\|#1\right\|}
\newcommand{\abs}[1]{\left\vert#1\right\vert}
\newcommand{\spr}[1]{\left\langle\,#1\,\right\rangle}

\newcommand{\Kl}[1]{\left\{#1\right\}}

\newcommand{\R}{\mathbb{R}} 

\newcommand{\N}{\mathbb{N}}

\def\div{\mfunc{div}} 

\newcommand{\xd}{x^\dagger}

\newcommand{\xad}{x_\alpha^\delta}

\newcommand{\xkd}{x_k^\delta}
\newcommand{\xk}{x_k}

\newcommand{\xkpd}{x_{k+1}^\delta}

\newcommand{\yd}{y^{\delta}}

\DeclareMathOperator*{\argmin}{argmin}

\newcommand{\ks}{{k_*}}
\newcommand{\xksd}{x_\ks^\delta}
\newcommand{\xt}{\tilde{x}}

\newcommand{\Lt}{L_2}

\newcommand{\kopt}{k_{\text{opt}}}
\newcommand{\Tad}{T_\alpha^\delta}

\newcommand{\as}{\alpha_*}

\newcommand{\kmax}{k_{\text{max}}}

\newcommand{\hr}{\psi}
\renewcommand{\div}{\operatorname{div}}

\newcommand{\Go}{{\small Good}}
\newcommand{\Ba}{{\small Bad}}
\newcommand{\Av}{{\small Average}}
\newcommand{\Ex}{{\small Excellent}}

\title{A numerical comparison of some heuristic stopping rules for nonlinear Landweber iteration}

\author{S.\ Hubmer\footnote{Johann Radon Institute Linz, Altenbergerstra{\ss}e 69, A-4040 Linz, Austria (simon.hubmer@ricam.oeaw.ac.at), corresponding author.} ,
E.\ Sherina\footnote{University of Vienna, Department of Mathematics, Oskar Morgenstern-Platz 1, 1090 Vienna, Austria (ekaterina.sherina@univie.ac.at)},
S.\ Kindermann\footnote{Johannes Kepler University Linz, Industrial Mathematics Institute, Altenbergerstra{\ss}e 69,A-4040 Linz, Austria (kindermann@indmath.uni-linz.ac.at)},
K.\ Raik\footnote{Johannes Kepler University Linz, Industrial Mathematics Institute, Altenbergerstra{\ss}e 69,A-4040 Linz, Austria (kemal.raik@indmath.uni-linz.ac.at)}}

\begin{document}

\maketitle

\begin{abstract}
The choice of a suitable regularization parameter is an important part of most regularization methods for inverse problems. In the absence of reliable estimates of the noise level, heuristic parameter choice rules can be used to accomplish this task. While they are already fairly well-understood and tested in the case of linear problems, not much is known about their behaviour for nonlinear problems and even less in the respective case of iterative regularization. Hence, in this paper, we numerically study the performance of some of these rules when used to determine a stopping index for Landweber iteration for various nonlinear inverse problems. These are chosen from different practically relevant fields such as integral equations, parameter estimation, and tomography.


\end{abstract}


\section{Introduction}\label{introduction}

In this paper, we consider nonlinear inverse problems of the form
    \begin{equation}\label{Fx=y}
	    F(x)=y \,, 
	\end{equation}
where $F:D(F)\subset X \to Y$ is a continuously Fr\'echet-differentiable nonlinear operator between real Hilbert spaces $X$ and $Y$. Furthermore, we assume that instead of exact data $y$ we are only given noisy data $\yd$ which satisfy
    \begin{equation}\label{cond_noise_delta}
        \norm{y-\yd} \leq \delta \,,
    \end{equation}
where $\delta$ denotes the noise level. Typically, inverse problems are  also ill-posed, which means that they may have no or even multiple solutions and in particular that a solution does not necessarily depend continuously on the data. This entails a number of difficulties, due to which one usually has to regularize the problem.

During the last decades, a large number of different regularization approaches have been developed; see for example \cite{Engl_Hanke_Neubauer_1996,Kaltenbacher_Neubauer_Scherzer_2008,Schuster_Kaltenbacher_Hofmann_Kazimierski_2012} and the references therein. Two of the most popular methods, which also serve as the bases for a wide variety of other regularization approaches, are \emph{Tikhonov regularization} \cite{Tikhonov_1963,Tikhonov_Glasko_1969} and \emph{Landweber iteration} \cite{Landweber_1951}. In its most basic form, Tikhonov regularization determines a stable approximation $\xad$  to the solution of \eqref{Fx=y} as the minimizer of the Tikhonov functional
    \begin{equation}\label{Tikhonov}
        \Tad(x) :=  \norm{F(x) - \yd}_Y^2 + \alpha \norm{x-x_0}_X^2 \,,
    \end{equation}
where $\alpha \geq 0$ is a regularization parameter and $x_0$ is an initial guess. In order to obtain convergence of $\xad$ to a solution $x_*$ of \eqref{Fx=y}, the regularization parameter $\alpha$ has to be suitably chosen. If the noise level $\delta$ from \eqref{cond_noise_delta} is known, then one can either use \emph{a-priori} parameter choice rules such as $\alpha \sim \delta$, or \emph{a-posteriori} rules such as the \emph{discrepancy principle}, which determines $\alpha$ as the solution of the nonlinear equation with some $\tau >1$
    \begin{equation}\label{discrepancy_principle_alpha}
        \norm{F(\xad)-\yd} = \tau \delta \,.  
    \end{equation}
Unfortunately, in many practical applications, estimates of the noise level $\delta$ are either unavailable or unreliable, which renders the above parameter choice rules impractical. Hence, a number of so-called \emph{heuristic parameter choice rules} have been developed over the years; see for example \cite{Leonov_1978,Leonov_1991,Tikhonov_Glasko_1965,Hanke_Raus_1996,Reginska_1996,Hansen_OLeary_1993,Wahba_1990} and the references therein. Most of them determine a regularization parameter $\as$ via
    \begin{equation}\label{heuristic_rule_alpha}
        \as \in \argmin_{\alpha \geq 0} \, \hr(\alpha,\yd) \,,
    \end{equation}
where $\hr : \R_0^+ \times Y \to \R \cup \{\infty\}$ is some lower semi-continuous functional. For example, the following popular choices in turn define the \emph{heuristic discrepancy (HD) principle} \cite{Hanke_Raus_1996}, the \emph{Hanke-Raus (HR) rule} \cite{Hanke_Raus_1996,Raus_Haemarik_2018}, the \emph{quasi-optimality (QO) rule} \cite{Tikhonov_Glasko_1969}, and the \emph{simple L (LS) rule} \cite{Kindermann_Raik_2020}:
	\begin{equation}\label{heuristic_rules_alpha}
	\begin{split}
	    \hr_{\text{HD}}(\alpha,y^\delta) &:= \frac{1}{\sqrt{\alpha}}\norm{F(\xad)-\yd}   \,,
	    \\
	    \hr_{\text{HR}}(\alpha,\yd) &:=\frac{1}{\alpha}\spr{ F(x^\delta_{\alpha,2})-\yd,F(\xad)-\yd}\,,
	    \\
	    \hr_{\text{QO}}(\alpha,\yd) &:= \norm{x^\delta_{\alpha,2}-\xad }\,,
	    \\
	    \hr_{\text{LS}}(\alpha,\yd)  &:=\spr{\xad-x^\delta_{\alpha,2},\xad} \,,
	\end{split}
	\end{equation}
where $x^\delta_{\alpha,2}$ denotes the so-called \emph{second Tikhonov iterate} \cite{Hanke_Groetsch_1998}, which is defined by
    \begin{equation*}
        x^\delta_{\alpha,2} := \argmin_{x\in X}\Kl{\norm{F(x)-\yd}_Y^2+\alpha\norm{x-\xad}_X^2 } \,.
    \end{equation*}
As is usually the case, rules of the form \eqref{heuristic_rule_alpha} are not restricted to Tikhonov regularization but can be used (possibly in a different form) in conjunction with various other regularization method as well. Due to their practical success in the treatment of linear inverse problems, heuristic parameter choice rules have received extended theoretical attention in recent years; see e.g.\ \cite{Kindermann_2011,Kindermann_2013,Kindermann_Neubauer_2008,Raus_Haemarik_2018,Haemarik_Palm_Raus_2011,Kindermann_Raik_Hamarik_Kangro_2019} and Section~\ref{sect_heuristics} for an overview.

A potential drawback of using Tikhonov regularization is the need to minimize the Tikhonov functional \eqref{Tikhonov}. Although for linear operators this reduces to solving a linear operator equation, the case of nonlinear operators is much more involved, since one typically has to use iterative optimization algorithms for the minimization of \eqref{Tikhonov}. Moreover, in order to use either the discrepancy principle \eqref{discrepancy_principle_alpha} or a heuristic rule of the form \eqref{heuristic_rule_alpha}, this minimization usually has to be done repeatedly for many different values of $\alpha$, which can render it infeasible for many practical applications.

Hence, a popular alternative in order to circumvent these issues is to directly use so-called \emph{iterative regularization methods}. As noted above, perhaps the most popular of these methods is Landweber iteration \cite{Engl_Hanke_Neubauer_1996,Kaltenbacher_Neubauer_Scherzer_2008}, which is defined by
	\begin{equation}\label{Landweber}
	    \xkpd = \xkd + \omega F'(\xkd)^*(\yd-F(\xkd)) \,,
	\end{equation}    
where $\omega > 0$ is a stepsize parameter. In order to obtain a convergent regularization method, Landweber iteration has to be combined with a suitable stopping rule such as the discrepancy principle, which now determines the stopping index $\ks$ by
    	\begin{equation}\label{discrepancy_principle}
	    \ks:=\min\Kl{k\in\N \, \vert \, \norm{F(\xkd)-\yd}\leq \tau\delta } \,,
	\end{equation}
for some parameter $\tau \geq 1$. Note that in contrast to \eqref{discrepancy_principle_alpha} for the choice of $\alpha$ in Tikhonov regularization, the discrepancy principle for Landweber iteration can be verified directly during the iteration, and does not require it to be run more than once. 

Now, analogously to \eqref{heuristic_rule_k}, most heuristic parameter choice (stopping) rules for Landweber iteration determine a stopping index $\ks$ via
    \begin{equation}\label{heuristic_rule_k}
        \ks \in \argmin_{k\in\N} \hr(k,\yd) \,,
    \end{equation}
with $\hr:\N\times Y\to\R \cup\{\infty\}$ again being some lower semi-continuous functional~\cite{Kindermann_2011,Kindermann_Neubauer_2008,Palm_2010,Kindermann_2013,Hanke_Raus_1996}. At least conceptually, the regularization parameter $\as$ in Tikhonov regularization and the stopping index $\ks$ for Landweber iteration play inversely proportional roles, i.e., $\alpha \sim 1/k$. Hence, the heuristic rules from \eqref{heuristic_rule_alpha} now correspond to
	\begin{equation}\label{heuristic_rules}
	\begin{split}
	    \hr_{\text{HD}}(k,y^\delta) &:=\sqrt{k}\norm{F(\xkd)-\yd}   \,,
	    \\
	    \hr_{\text{HR}}(k,\yd) &:=k \spr{\yd-F(x^\delta_{2k}), \yd-F(\xkd)} \,,
	    \\
	    \hr_{\text{QO}}(k,\yd) &:= \norm{x^\delta_{2k} - \xkd} \,,
	    \\
	    \hr_{\text{LS}}(k,\yd)  &:= \spr{\xkd, x^\delta_{2k} - \xkd } \,.
	\end{split}
	\end{equation}
Note that the analogue of the second Tikhonov iteration $x^\delta_{\alpha,2}$ in the Landweber method corresponds to applying $k$ steps of a Landweber iteration with initial guess $\xkd$, which is easily seen to be identical to simply doubling the iteration number, i.e., calculating $x^\delta_{2k}$. 

Similarly to the case of the discrepancy principle \eqref{discrepancy_principle}, these functionals can be evaluated during a single run of Landweber iteration. Hence, in particular in the nonlinear case, the combination of heuristic parameter choice rules together with Landweber iteration (or other iterative regularization methods) suggests itself.

Heuristic parameter choice rules are already fairly well-understood in the case of linear problems. Despite a number of obstacles such as \emph{Bakushinskii's veto} \cite{Bakushinskiy_1985}, many theoretical results on both the convergence and other aspects of these rules are already available; see, e.g.,~\cite{Kindermann_2011} and the references therein. Furthermore, numerical tests have been carried out for various different linear test problems \cite{Bauer_Lukas_2011,Haemarik_Palm_Raus_2011,Palm_2010}. In contrast, for the case of nonlinear problems, not much is known with respect to convergence theory 
nor is there a unifying framework for heuristic rules, and, furthermore, there are no numerical performance studies for heuristic rules. 
 
Hence, the main motivation of this article is to close this gap. In particular, we consider the behaviour of the four rules given in \eqref{heuristic_rule_k}, i.e., the heuristic discrepancy principle, the Hanke-Raus rule, the quasi-optimality rule, and the simple L rule, using Landweber iteration and for different nonlinear test problems from practically relevant fields such as integral equations, parameter estimation, and tomography. The aim of this paper is two-fold: On the one hand, we want to demonstrate that heuristic parameter choice rules can indeed be used successfully not only for linear but also for nonlinear inverse problems. On the other hand, we want to provide some useful insight into potential difficulties and pitfalls which one might encounter, and what can be done about them. Since we also compare the heuristics with the discrepancy principle, the results in this article additionally provide a numerical study of the performance  of the latter in the nonlinear case.

The outline of this paper is as follows: in Section~\ref{sect_background}, we provide some general background on heuristic parameter choice rules and nonlinear Landweber iteration. In Section~\ref{sect_test_problems}, we then introduce the various problems on which we want to test the different heuristic parameter choice rules. The corresponding results are then presented in Section~\ref{sect_numerical_results}, which is followed by a short conclusion in Section~\ref{sect_conclusion}.

\section{Landweber iteration and stopping rules}\label{sect_background}

In this section, we recall some basic results on Landweber iteration and heuristic parameter choice rules. Since in this paper we are mainly interested in a numerical comparison of the different rules, here we only provide a general overview of some of the main results from the literature, focusing in particular on those aspects which are relevant to understand the numerical results presented below.

\subsection{Landweber iteration for nonlinear problems}\label{sect_Landweber}

One of the main differences to the linear case is that for nonlinear problems only local convergence can be established. For this, one needs to impose certain restrictions on the nonlinearity of the operator $F$, such as the \emph{tangential cone condition} \cite{Hanke_Neubauer_Scherzer_1995,Kaltenbacher_Neubauer_Scherzer_2008}, which is given by
    \begin{equation}\label{tangentialconecondition}
	    \norm{F(x)-F(\xt)-F'(x)(x-\xt)} \leq \eta \norm{F(x)-F(\xt)}
	    \,,
	    \qquad
	    \eta< 1/2 \,.
	\end{equation}
This condition has to hold locally in a neighbourhood of a solution $x^*$ of the problem, and the initial guess $x_0$ of the iteration has to be contained inside it. Furthermore, the parameter $\tau$ in the discrepancy principle \eqref{discrepancy_principle} has to satisfy
	\begin{equation}\label{cond_tau}
	    \tau > 2\frac{1+\eta}{1-2\eta} \geq 2 \,,
	\end{equation}
where the factor $2$ can be slightly improved by an expression depending on $\eta$ which tends to $1$ as $\eta \to 0$, thereby recovering the optimal bound in the linear case \cite{Hanke_2014}. If in addition the stepsize $\omega$ is chosen small enough such that locally there holds
	\begin{equation}\label{cond_omega}
	    \omega\norm{F'(x)}^2 \leq 1 \,,
	\end{equation}
then Landweber iteration combined with the discrepancy principle \eqref{discrepancy_principle} converges to  $x^*$ as the noise level $\delta$ goes to $0$ \cite{Hanke_Neubauer_Scherzer_1995,Kaltenbacher_Neubauer_Scherzer_2008}. Furthermore, convergence to the minimum-norm solution $\xd$ can be established given that in a sufficiently large neighbourhood there holds $N(F'(\xd))\subset N(F'(x))$. In order to prove convergence rates, in addition to source conditions further restrictions on the nonlinearity of $F$ are necessary \cite{Engl_Hanke_Neubauer_1996,Kaltenbacher_Neubauer_Scherzer_2008}.

Even though the tangential cone condition \eqref{tangentialconecondition} holds for a number of different applications (see e.g.~\cite{Kaltenbacher_Nguyen_Scherzer_2019} and the references therein), and even though attempts have been made to replace it by more general conditions \cite{Kindermann_2017}, these can still be difficult to prove for specific applications (cf.,~e.g., \cite{Kindermann_2022} for an analysis for the EIT problem and \cite{Hubmer_Sherina_Neubauer_Scherzer_2018} for an analysis of a parameter estimation problem in linear elastography). Furthermore, even if the tangential cone condition can be proven, the exact value of $\eta$ typically remains unknown. Since this also renders condition \eqref{cond_noise_delta} impractical, the parameter $\tau$ in the discrepancy principle then has to be chosen manually; popular choices include $\tau = 1.1$ or $\tau = 2$. These work well in many situations, but are also known to fail in others (compare with Section~\ref{sect_numerical_results} below). In any case, this shows that for practical applications involving nonlinear operators, informed ``heuristic'' parameter choices remain necessary even if the noise level $\delta$ is known.

\subsection{Heuristic stopping rules}\label{sect_heuristics}

As mentioned in the introduction, in many practical situations one does not have knowledge of the noise level. Thus, applying it with unreliable estimates of the noise level is rarely fruitful. The remedy is the use of \emph{heuristic} (aka \emph{data-driven} or \emph{error-free}) rules, where the iteration is terminated at $\ks:=k(\yd)$, which depends only on the measured data and not the noise level. 

The analysis of the so-called heuristic rules for Tikhonov regularisation has been examined extensively for linear problems \cite{Kindermann_2011,Kindermann_Neubauer_2008,Kindermann_2013} but less so beyond the linear case. Some results for convex problems and nonlinear problems in Banach spaces can be found in \cite{Jin_Lorenz_2010,Jin_2016_02,Jin2017}. The numerical performance of these rules in the linear case has also been studied \cite{Bauer_Lukas_2011,Haemarik_Palm_Raus_2011,Palm_2010} and for convex Tikhonov regularization in \cite{Kindermann_Mutimbu_Resmerita_2014}. However, for Landweber iteration for nonlinear problems, neither an analysis has been given nor has the numerical performance of heuristic stopping rules been investigated so far.  

Let us briefly illustrate the rationale behind the heuristic minimization-based rules \eqref{heuristic_rule_k} (or \eqref{heuristic_rule_alpha}). It is well-known that the total error between the regularized solution $\xkd$ and the  exact solution $\xd$ can be split into approximation and stability error:
    \begin{equation}\label{eq_error_split}
        \norm{\xkd- \xd} \leq \norm{\xk - \xd} + \norm{\xk-\xkd} \,. 
    \end{equation}
Here, $\xk$ denotes the regularized solution when exact data ($\delta = 0$) would be given. An ideal optimal parameter choice would be one that minimizes over $k$ the total error or the upper bound in \eqref{eq_error_split}. However, this ``oracle'' parameter choice is not possible, as in \eqref{eq_error_split} only the element $\xkd$ is at hand, and neither the exact solution nor the exact data are available. The idea of minimization-based heuristic rules is to construct a computable  functional $\hr$ that estimates the total error sufficiently well, i.e., $\hr(k,\yd) \sim \norm{\xkd- \xd}$, such that a minimization over $k$ is expected to yields a reasonable parameter choice of $\ks$ as well. In analogy to 
the approximation/stability split, we may just as well use a similar splitting for $\hr$, i.e., 
$\hr(k,\yd) \leq \hr_a(k,\yd) + \hr_d(k,\yd)$ and search for conditions such that the respective parts estimate the approximation and stability error separately:
    \begin{align} 
        \hr_a(k,\yd) \sim \norm{\xk - \xd} \,, \label{eq_psiapprox} 
        \\ 
        \hr_d(k,\yd) \sim \norm{\xk-\xkd } \,.  \label{eq_psistab} 
    \end{align} 

As briefly mentioned earlier, the pitfall for heuristic stopping rules manifests itself in the form of the so-called \emph{Bakushinskii veto}, the consequence of which is that a heuristic stopping rule cannot yield a convergent regularisation scheme in the \emph{worst case} scenario \cite{Bakushinskiy_1985}, i.e., for all possible noise elements $y-\yd$. A direct consequence is that there \emph{cannot} exist a $\hr$ with the error-estimating capabilites as mentioned above, in particular, such that \eqref{eq_psistab} holds!

However, there is a way to overcome the negative result of Bakushinskii by restricting the class of permissible noise elements $y-\yd$. In this way, one can prove convergence in a \emph{restricted noise case} scenario. (Of course, for this approach to be meaningful, the restrictions should be such that ``realistic" noise is always permitted.) Some noise restriction were used, e.g., in \cite{Glasko_Kriskin_1984,Hanke_Raus_1996}, although the restrictions there were implicit and very hard to interpret. (For instance, in \cite{Glasko_Kriskin_1984}, essentially, condition \eqref{eq_psistab} was \emph{postulated} rather than derived from more lucid conditions.)

A major step towards an understanding of heuristic rules was made in \cite{Kindermann_Neubauer_2008,Kindermann_2011}, when a full convergence analysis in the linear case with explicit and interpretable restrictions was given. These conditions, which were proven to imply
\eqref{eq_psistab}, take the form of a  \emph{Muckenhoupt}-type inequality: Let $(\sigma_i,u_i,v_i)$ be the singular system of the forward operator. Then, for $p \in \{1,2\}$, the $p$-Muckenhoupt-inquality holds if there is a constant $C$ and a $t_0$ such that for all admissible noise $y-\yd$ it holds that 
    \begin{equation}\label{eq_mc1}  
        \sum_{\sigma_i^2 \geq t}  
        \frac{t}{\sigma_i^2} \abs{(y-\yd,v_i)}^2  
        \leq C  
        \sum_{\sigma_i^2 < t}  \left(\tfrac{\sigma_i^2}{t}\right)^{p-1}  \abs{(y-\yd,v_i)}^2 \qquad \forall t \in (0,t_0).
    \end{equation} 
It has been shown in \cite{Kindermann_Neubauer_2008,Kindermann_2011,Kindermann_Raik_2020} that for most of the classical regularization schemes and for the four above mentioned rules  this condition suffices for being able to estimate the stability error  
\eqref{eq_psistab}, and thus convergence can be proven. The inequality \eqref{eq_mc1} has the interpretation 
of an  ``irregularity'' condition for the noise vector $y-\yd$; by postulating \eqref{eq_mc1}, the noise must be distinguishable from smooth data error (which never satisfies \eqref{eq_mc1}). However, this anyway agrees with the common idea of noise. 

\begin{remark}\label{rem_heuremark}
The above heuristic rules require slightly different Muckenhoupt conditions
which leads to two groups of rules: For HD and HR, \eqref{eq_mc1} with $p = 1$ suffices, while for QO and LS, the condition with $p =2$ (which is a slightly stronger requirement) has to be postulated~\cite{Kindermann_Neubauer_2008,Kindermann_2011,Kindermann_Raik_2020}. Thus the former might be successful even if the later fail. However, it has to be kept in mind that the error analysis shows that, as long as they 
can be successfully applied, QO and LS in general lead to smaller errors. 
\end{remark}

Indeed, it has been shown in \cite{Kindermann_Neubauer_2008} that the above mentioned Muckenhoupt inequality is satisfied in typical situations, and in \cite{Kindermann_Pereverzyev_Philipenko_2018} it was shown in a stochastic setting that it is also satisfied for coloured Gaussian noise almost surely for many cases. Below, we discuss cases where the Muckenhoupt inequality might not be satisfied and when heuristic rules may fail.

The above mentioned noise restrictions are heavily rooted in the linear theory and in particular make use of spectral theory and the functional calculus of operators. In the case of a nonlinear operator, we are no longer afforded the luxury of having these tools available. Some alternative noise conditions in the nonlinear case have been established in~\cite{Jin_Lorenz_2010} for convex variational Tikhonov regularisation, in \cite{Zhang_Jin_2018} for Bregman iteration in Banach spaces, or in~\cite{Liu_Real_Lu_Jia_Jin_2020} for general variational regularisation. However, as of yet, these conditions could not be  deciphered into a palatable explanation as to when a rule will work or otherwise. An attempt was made in \cite{Kindermann_Raik_2019_02} to formulate an analogous Muckenhoupt-type inequality for convex Tikhonov regularisation in a somewhat restrictive setting of a diagonal operator over $\ell^p$ spaces with $p\in(1,\infty)$. 

However, to the knowledge of the authors, neither a convergence analysis  nor a numerical investigation of heuristic rules for nonlinear Landweber iteration seems to be available in the literature. In light of all of this, we thus have great incentive to investigate heuristic stopping rules for nonlinear Landweber iteration numerically and to compare them with the more tried and tested a-posteriori rules. 

The Muckenhoupt inequality covers the convergence theory. When it comes to the practical capabilities of heuristic rules, also convergence rates are important. For this, efficient estimates of the approximation error as in \eqref{eq_psiapprox} are vital, and in the linear case, sufficient conditions for this have been established, this time in form of conditions for the exact solution $\xd$. For instance, using the singular system $(\sigma_i,u_i,v_i)$, the following {\em regularity condition} (for $p = 1$ or $p=2$, depending on 
the rule) 
    \begin{equation}\label{eq_reg_con} 
        \sum_{\sigma_i^2 \leq t } \abs{(\xd ,u_i)}^2 \leq 
        C  \sum_{\sigma_i^2 > t } \left(\tfrac{\sigma_i^2}{t}\right)^{p-1}\!
        r(\sigma_i^2,\alpha)^2
        \abs{(\xd ,u_i)}^2, \qquad \forall t \in (0,t_0) \,,
    \end{equation}
is sufficient for \eqref{eq_psiapprox} and, together with smoothness conditions, yields (optimal-order) convergence rates in many situations \cite{Kindermann_2011}. Here, $r(\lambda,\alpha)$ is the residual spectral filter function of the regularization method, i.e., $r(\lambda,\alpha) = \frac{\alpha}{\alpha +\lambda}$ for Tikhonov regularization and $r(\lambda,\alpha) = (1-\lambda)^\frac{1}{\alpha}$ with $\alpha = k^{-1}$ for Landweber iteration. Note that the regularity condition depends strongly on the regularization method via $r(\lambda,\alpha)$ in contrast to the Muckenhoupt condition. The rough interpretation of \eqref{eq_reg_con} is that the exact solution has coefficients $(\xd ,u_i)$ that do not deviate too much from a given decay (that is encoded in a smoothness condition).

\subsection{Challenges and practical issues for heuristic rules}\label{subsect_challenges}

Next, let us point out possible sources of failure for heuristic rules and some  peculiarities for the case of nonlinear Landweber regularization. 

\begin{itemize} 
    \setlength{\itemsep}{5pt}
    \item \emph{Failure of Muckenhopt condition.} A general problem for heuristic rules, both in the linear and nonlinear case, is when the Muckenhoupt condition \eqref{eq_mc1} is not satisfied. This can happen for standard noise for super-exponential ill-posed problems, for instance, for the backward heat equation (see \cite{Kindermann_2011,Kindermann_2013}). Less obvious but practically important is the case that the Muckenhoupt inequality might also fail to hold for standard noise if the problem is \emph{nearly well-posed}, i.e., when the singular values decay quite slowly (e.g., as $\sigma_i \sim i^{-\beta}$ with, say, $\beta < 1$). In this case, exact data can be quite irregular as the operator is only little smoothing, and it is hard to distinguish between exact data and noise, and this is indeed a relevant possible source of failure. 
    \item \emph{The spurious first local minimum for Landweber iteration.} Recall that the effective performance of heuristic rules depend also on efficent estimates of the approximation errors and  in particular on the regularity condition to hold. As pointed out before, this strongly depends on the regularization method. We have extensive numerical evidence that for linear and nonlinear Landweber iteration, the approximation error (i.e., \eqref{eq_psiapprox}) is often only badly estimated by $\hr(k,\yd)$ for the first few iterations (i.e., $k$ small) and that \eqref{eq_reg_con} holds only with a bad constant for large $t = k^{-1}$. In practical computations, this has the consequence that  $\hr(k,\yd)$ typically has an outstanding local minimum for small $k$. However, this local minimum is rarely the global minimum  which usually  appears much later for larger $k$, and inexperienced users are often tempted to mistakenly take this local minimum for the global one to save having to compute later iterations. This happens quite often in the linear and the nonlinear case for Landweber iteration, but a similar problem for Tikhonov regularization is rarely observed. The deeper reason for this discrepancy is the different shape of the residual filter function for both method, which makes the regularity condition \eqref{eq_reg_con} more restrictive for Landweber iteration and for large $t$. 
    \item \emph{Discretization cut-off.} It is known that, due to discretization, the theoretically global minimum  of $\hr(k,\yd)$ for finite-dimensional problems is at $k= \infty$, which does not provide a correct stopping index. Thus, in practical computations, we have to restrict the search space by fixing an upper bound for $k$ (or, for continuous regularization method, by a lower bound for  $\alpha$). Some rules how to do this together with an accompanying analysis in the linear case are given in \cite{Kindermann_2013}; however,  for nonlinear  problems no such investigation exists. This issue is  relevant for very small noise levels or for coarse discretizations, and in practice one takes a pragmatic approach  and assumes a reasonable upper bound for the iteration index and looks for interior minima rather than global minima at the boundary of the search space. 
    \item \emph{Only local convergence in nonlinear case.} The established convergence theory in the nonlinear case is a local one: one can only prove convergence when the initial guess is sufficiently close to the exact solution $\xd$, and in the case where noise is present, the iteration usually diverges out of the neighborhood of $\xd$ as $k\to \infty$. In particular, it is possible that $\xkd$ ``falls'' out of the domain of the forward operator. As a consequence, it might happen that the functionals $\hr(k,\yd)$ in \eqref{heuristic_rule_k} are not defined for very large $k$. By definition, however, one would have to compute a minimizer over all $k$, which  is then not practically possible. (This is different to Tikhonov regularization, whose solution is always well-defined for any $\alpha$). In practice, as a remedy, one would introduce an upper limit for the number of iterations up to which the functional $\hr(k,\yd)$ is computed. Additionally, one could monitor the distance to the intial guess $\norm{\xkd-x_0}$ and terminate if this becomes too large.  
\end{itemize} 

In this section, have stated some practical aspects of heuristic rules; a deeper mathematical analysis (especially in the nonlinear case) is outside the scope of this article. For further aspects on heuristic stopping rules both from a theoretical and practical viewpoint, we refer the reader to \cite{Raik_20} and the references therein.

\section{Test problems}\label{sect_test_problems}
	
In this section, we introduce a number of test problems on which we evaluate the performance of the heuristic stopping rules described above. These nonlinear inverse problems belong to a variety of different problem classes, including integral equations, tomography, and parameter estimation. For each of them, we shortly review their background and describe their precise mathematical setting and relevant theoretical results below.

\subsection{Nonlinear Hammerstein operator}\label{subsect_Hammerstein}
	
A commonly used nonlinear inverse problem \cite{Hanke_Neubauer_Scherzer_1995, Neubauer_2000,Neubauer_2016, Neubauer_2017_2, Hubmer_Ramlau_2017, Hubmer_Ramlau_2018} for testing, in particular, the behaviour of iterative regularization methods is based on so-called nonlinear Hammerstein operators of the form
    \begin{equation*}
        F \, : \, H^1[0,1] \to \Lt[0,1] \,, 
        \qquad
        F(x)(s) := \int_0^1 k(s,t) \gamma(x(t)) \,dt \,,
    \end{equation*}
with some given function $\gamma: \R\to \R$.    
Here, we look at a special instance of this operator, namely
    \begin{equation}\label{def_Hammerstein}
        F(x)(s) := \int_0^s x(t)^3 \, dt \,,
    \end{equation}
for which the tangential cone condition \eqref{tangentialconecondition} holds locally around a solution $x^\dagger$, given that it is bounded away from zero (see e.g.\ \cite{Neubauer_2017_2}). Furthermore, the Fr\'echet derivative and its adjoint, which are required for the implementation of Landweber iteration, can be computed explicitly.

\subsection{Diffusion-Coefficient estimation}\label{subsect_diffusion}
	
Another classic test problem \cite{Kaltenbacher_Neubauer_Scherzer_2008} in inverse problems is the estimation of the diffusion coefficient $a$ in the partial differential equation
    \begin{equation*}
        - \div(a \nabla u) = f \,,
    \end{equation*}
from measurements of $u$, and given knowledge of the source-term $f$ and (Dirichlet) boundary conditions on $u$. For this test problem, we focus on the one-dimensional version
    \begin{equation}\label{eq_diffusion_PDE}
    \begin{split}
        -(a(s)u(s)_s)_s = f(s) &\,, \qquad s \in (0,1) \,,
        \\
        u(0) = u(1) = 0 &\,,
    \end{split}
    \end{equation}
which leads to an inverse problem of the form \eqref{Fx=y} with the nonlinear operator
    \begin{equation}\label{def_diffusion}
    \begin{split}
        F \, : \, D(F) := \{ a \in H^1[0,1] \, : \, a(s) \geq \underline{a} > 0  \} \quad &\to \quad L^2[0,1] \,,
        \\
        a \quad &\mapsto \quad F(a) := u(a) \,,
    \end{split}        
    \end{equation}
where $u(a)$ is the solution of \eqref{eq_diffusion_PDE} above. The computation of the Fr\'echet derivative and its adjoint of $F$ now requires solving PDEs of the form \eqref{eq_diffusion_PDE}. Furthermore, it was shown (see e.g.\ \cite{Kaltenbacher_Neubauer_Scherzer_2008}) that the tangential cone condition \eqref{tangentialconecondition} holds locally around a solution $a^\dagger \geq c > 0$.

\subsection{Acousto-Electical Tomography}\label{subsect_AET}
	
Another PDE parameter estimation problem, this time from the field of tomography, is the hybrid imaging modality of acousto-electrical tomography (AET) \cite{Ammari_Bonnetier_Capdeboscq_Tanter_Fink_2008, Kuchment_Kunyansky_2010,Gebauer_Scherzer_2008, Zhang_Wang_2004}. Based on a modulation of electrical impedance tomography (EIT) by ultrasound waves, AET aims at reconstructing the spatially varying electrical conductivity distribution inside an object from electrostatic measurements of voltages and the corresponding current fluxes on its surface. Compared for example to EIT, reconstructions of high contrast and high resolution may be obtained. Mathematically, the problem amounts to reconstructing the spatially varying conductivity $\sigma$ from measurements of the power densities 
    \begin{equation*}
        E_j(\sigma) := \sigma \abs{ u_j(\sigma) }^2 \,,
    \end{equation*}
where the interior voltage potentials $u_j(\sigma)$ are the solution of the elliptic PDEs 
    \begin{equation}\label{eq_AET_PDE}
    \begin{split}
        \text{div}(\sigma \nabla u_j) &= 0 \,,
        \quad \text{in}\, \Omega \,,
        \\
        (\sigma \nabla u_j) \cdot \Vec{n} \vert_{\partial\Omega} &= g_j \,,
    \end{split}
    \end{equation}
where $\Omega \subset \R^N$, $N=2,3$ is a bounded and smooth domain, and $g_j$ models the current flux on the boundary $\partial \Omega$ in the outward unit normal direction $\Vec{n}$. Once again, this problem can be restated as an operator equation of the form \eqref{Fx=y} with a Fr\'echet differentiable nonlinear operator. Its Fr\'echet derivative and the adjoint thereof can for example be found in \cite{Hubmer_Knudsen_Li_Sherina_2018}, and their evaluation again require the solution of PDEs of the form \eqref{eq_AET_PDE} for different right-hand sides. Note that it is not known whether the tangential cone condition holds for AET (or EIT). Furthermore, it is in general not possible to uniquely determine the conductivity $\sigma$ from a single power density measurement $E_j(\sigma)$ \cite{Bal_2013,Isakov_2006}. In addition, if $g_j = 0$ on some part $\Gamma \subset \partial \Omega$ of the boundary, then the problem becomes severely ill-posed. On the other hand, if $g_j \neq 0$ almost everywhere on $\partial \Omega$ and given a sufficient amount (depending on the dimension $N$) of ``different'' power density measurements $E_j(\sigma)$, the conductivity $\sigma$ can be uniquely reconstructed \cite{Capdeboscq_Fehrenbach_Gournay_Kavian_2009, Bal_Bonnetier_Monard_Triki_2013, Monard_Bal_2012, Alberti_Capdeboscq_2018}. In this case, the problem behaves numerically close to well-posed, which is reflected in the behaviour of the heuristic parameter choice rules; cf.~Section~\ref{subsect_challenges}.

\subsection{SPECT}\label{subsect_SPECT}

Next, we look at Single Photon Emission Computed Tomography (SPECT), which is another example from the large field of tomography \cite{Natterer_2001, Dicken_1998, Dicken_1999, Ramlau_2003, Ramlau_Teschke_2006}. In this medical imaging problem, on aims at reconstructing the radioactive distribution $f$ (activity function) and the attenuation map $\mu$, which is related to the density of different tissues, from radiation measurements outside the examined body. The connection between these quantities is typically modelled by the attenuated Radon Transform (ART), which is given by \cite{Natterer_2001}: 	
    \begin{equation}\label{def_SPECT}
        F(f,\mu)(s,\omega) := \int_{\R} f(s \omega^\perp + t \omega) \exp\left(-\int_{t}^{\infty} \mu(s \omega^\perp + r \omega) \, dr\right) \, dt \,,
    \end{equation} 
where $s \in \R$ and $\omega \in S^1$. With this, one again arrives at a problem of the form \eqref{Fx=y}, where $y$ then is the measured singoram. The well-definedness and differentiability of the operator $F$ with respect to suitable Sobolev spaces has been studied in detail in \cite{Dicken_1998, Dicken_1999}. However, it is still unknown whether also the tangential cone condition \eqref{tangentialconecondition} holds for this problem.

\subsection{Auto-Convolution}\label{subsect_conv}
	
As a final test example, we consider the problem of (de-)auto-convolution \cite{Buerger_Hofmann_2015,Fleischer_Hofmann_1996, Gorenflo_Hofmann_1994,Ramlau_2003}. Among the many inverse problems based on integral operators, auto-convolution is particularly interesting due to its importance in laser optics \cite{Anzengruber_Buerger_Hofmann_Steinmeyer_2016,Birkholz_Steinmeyer_Koke_Gerth_Buerger_Hofmann_2015, Gerth_2014}. Mathematically, it amounts to solving an operator equation of the form \eqref{Fx=y} with the operator
    \begin{equation}\label{def_conv}
        F \, : \, \Lt[0,1] \to \Lt[0,1] \,,
        \qquad
        F(x)(s) := (x \ast x) (s) := \int_0^1 x(s-t)x(t) \, dt \,,
    \end{equation}
where the functions on $\Lt[0,1]$ are interpreted as 1-periodic functions on $\R$. While deriving the Fr\'echet differentiability and its adjoint of $F$ is straightforward, it is not known whether the tangential cone condition \eqref{tangentialconecondition} holds. However, for small enough noise levels $\delta$, the residual functional is locally convex around the exact solution $\xd$, given that it only has finitely many non-zero Fourier coefficients \cite{Hubmer_Ramlau_2018}.

\section{Numerical Results}\label{sect_numerical_results}

In this section, we present the results of using the four heuristic parameter choice rules defined in \eqref{heuristic_rules} to determine a stopping index for Landweber iteration applied to the different nonlinear test problems introduced in Section~\ref{sect_test_problems}. 

For each of these problems, we started from a known solution $\xd$ in order to define the exact right-hand side $y$. Random noise corresponing to different noise levels $\delta$ was added to $y$ in order to create noisy data $\yd$, and a suitable stepsize $\omega$ for Landweber iteration was computed via \eqref{cond_omega} based on numerical estimates of $\norm{F'(\xd)}$. Afterwards, we ran Landweber iteration for a predefined number of iterations $\kmax$, which was chosen manually for each problem via a visual inspection of the error, residual, and heuristic functionals, such that all important features of the parameter choice rules were captured for this comparison. 

Following each application of Landweber iteration, we computed the values of the heuristic functionals $\hr$, as well as their corresponding minimizers $\ks$. As noted in Section~\ref{sect_background}, the functional values corresponding to the first few iterations have to be discarded in the search for the minimizers due to the spurious
first local minimum (tacitly assuming that the noise-level is small enough such 
that a good stopping index appears later). 
For each of the different heuristic rules we then computed the resulting absolute error 
	\begin{equation}\label{abs_error}
	    \norm{\xksd -  \xd }\,,
	\end{equation}
and for comparison, for each problem we also computed the \emph{optimal stopping index}
	\begin{equation}\label{def_kopt}
	    \kopt := \argmin_{k\in\N} \norm{\xkd-\xd} \,,
	\end{equation}
together with the corresponding optimal absolute error. Furthermore, we also computed the stopping index $k_\text{DP}$ determined by the discrepancy principle \eqref{discrepancy_principle}, which can also be interpreted as the ``first'' minimizer of the functional
    \begin{equation}\label{def_hrDP}
        \hr_\text{DP}(k) := \abs{ \norm{F(\xkd) - \yd} - \tau \delta} \,.
    \end{equation}
As noted in Section~\ref{sect_Landweber}, since the exact value of $\eta$ in \eqref{tangentialconecondition} is unknown for our test problems, a suitable value for $\tau$ has to be chosen manually. Depending on the problem, we used either one of the popular choices $\tau = 1.1$ or $\tau = 2$, although as we are going to see below, these are not necessarily the ``optimal'' ones. In any case, the corresponding results are useful reference points to the performance of the different heuristic parameter choice rules.     

Concerning the discretization and implementation of each of the numerical test problems, we refer to the subsequent sections and the references mentioned therein. All computations were carried out in Matlab on a notebook computer with an Intel(R) Core(TM) i7-85650 processor with 1.80GHz (8 cores) and 16 GB RAM, except for the acousto-electrical tomography problem, which was carried out in Python using the FEniCS library \cite{Alnaes_Blechta_2015} on a notebook computer with an Intel(R) Core(TM) i7-4810MQ processor with 2.80GHz (8 cores) with 15.3 GB RAM.

\subsection{Nonlinear Hammerstein operator}

\begin{figure}[ht!]
    \centering
    \includegraphics[width=\textwidth, trim = {6.5cm 1.5cm 6cm 2cm}, clip = true]{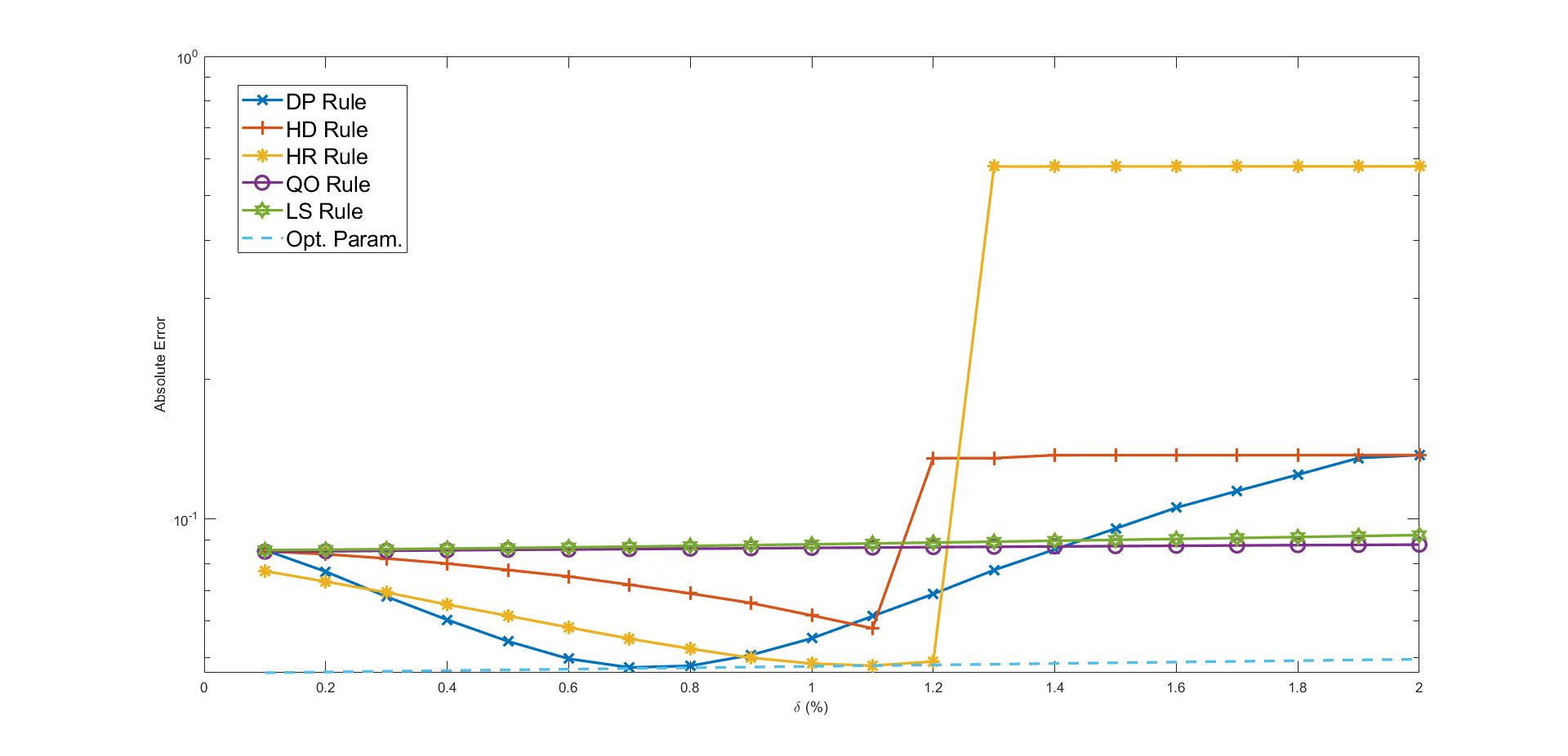}
    \caption{Numerical results for the nonlinear Hammerstein problem introduced in Section~\ref{subsect_Hammerstein}: Absolute error \eqref{abs_error} at the stopping indices $\ks$ determined by the discrepancy principle \eqref{discrepancy_principle}, the different heuristic parameter choice rules \eqref{heuristic_rules}, and the optimal stopping index $\kopt$ defined in \eqref{def_kopt}, each for different relative noise levels $\delta$.}
    \label{fig_Hammerstein_results}
\end{figure}

\begin{figure}[ht!]
    \centering
    \includegraphics[width=\textwidth, trim = {6.5cm 1.5cm 6cm 2cm}, clip = true]{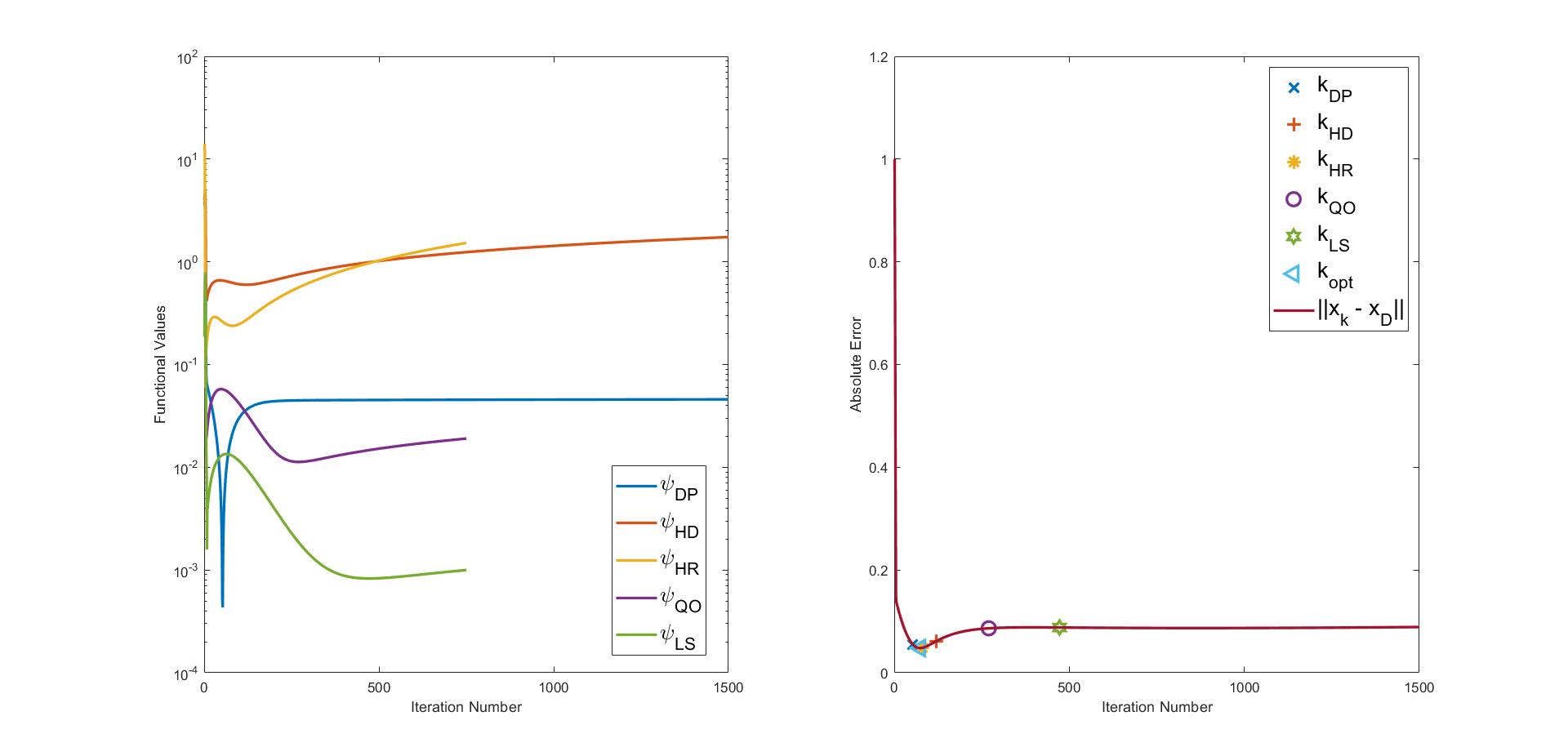}
    \caption{Numerical results for the nonlinear Hammerstein problem introduced in Section~\ref{subsect_Hammerstein}: Heuristic functionals \eqref{heuristic_rules} and discrepancy functional \eqref{def_hrDP} for $\delta = 1\%$ relative noise (left). Corresponding evolution of the absolute error $\norm{\xkd - \xd}$ with marked points indicating the stopping indices chosen by the different rules (right).}
    \label{fig_Hammerstein_functionals}
\end{figure}

First, we consider the nonlinear Hammerstein problem introduced in Section~\ref{subsect_Hammerstein}. In order to discretize this problem, the interval $[0,1]$ is subdivided into $128$ subintervals, and the operator $F$ itself is discretized as described in~\cite{Neubauer_2000,Neubauer_2017_2}; cf.\ also \cite{Hubmer_Ramlau_2017,Hubmer_2015}. For the exact solution we choose $\xd(s) = 2+(s-0.5)/10$ and compute the corresponding data $y$ by the  application of the Hammerstein operator \eqref{def_Hammerstein}. For the initial guess we choose $x_0(s) = 1$, and in the discrepancy principle \eqref{discrepancy_principle} we use $\tau = 2$. The absolute error \eqref{abs_error} corresponding to different parameter choice rules and noise levels $\delta$ from $0.1\%$ to $2\%$ is depicted in Figure~\ref{fig_Hammerstein_results}. Typical plots of the heuristic functionals $\hr$ as well as the evolution of the absolute error over the iteration are depicted in Figure~\ref{fig_Hammerstein_functionals}. There, the marked points denote the corresponding stopping indices selected via the different rules. 

As can be seen from the left plot in Figure~\ref{fig_Hammerstein_functionals}, the heuristic functionals $\hr$ generally exhibit the same shape as expected from the theoretical considerations discussed above. For example, each of the functionals exhibits a spurious ``first'' local minimum within the first few iterations, as already discussed in Section~\ref{subsect_challenges}. Apart from this, each functional $\hr$ has a well-defined minimum reasonably close to the stopping index $k_{\text{DP}}$ determined by the discrepancy principle. However, for larger noise levels, this minimum vanishes for the  the HD and the HR rule, which is reflected in Figure~\ref{fig_Hammerstein_results} by their unsatisfactory constant absolute error (the rules select the spurious minimum in this case). In contrast, the QO and LS rule keep their general shape for all noise levels, and thus produce stable stopping indices, which are typically larger than those determined by the discrepancy principle. Since the evolution of the absolute error depicted in the right plot in Figure~\ref{fig_Hammerstein_functionals} flattens for larger iteration numbers (an effect of the discretization), the error curves for the QO and the LS rules remain rather constant on the logarithmic scale. Curiously, and contrary to theoretical expectations, the absolute error curve for the discrepancy principle in Figure~\ref{fig_Hammerstein_results} exhibits a parabola shape. This indicates that on the one hand, the chosen value of $\tau$ is too small, while on the other hand the discretization might be too coarse in comparison with the small noise levels. However, in practice the discretization is often fixed by practical limitations, while a proper value of $\tau$ satisfying \eqref{cond_tau} is typically impossible to determine, or unreasonably large. Hence, this first test already indicates the usefulness of heuristic parameter choice rules (expecially the QO and LS rule) in comparison to the discrepancy principle, and shows some typical limitations which we now investigate further in the remaining test problems.

\subsection{Diffusion-Coefficient estimation}

\begin{figure}[ht!]
    \centering
    \includegraphics[width=\textwidth, trim = {6.5cm 1.5cm 6cm 2cm}, clip = true]{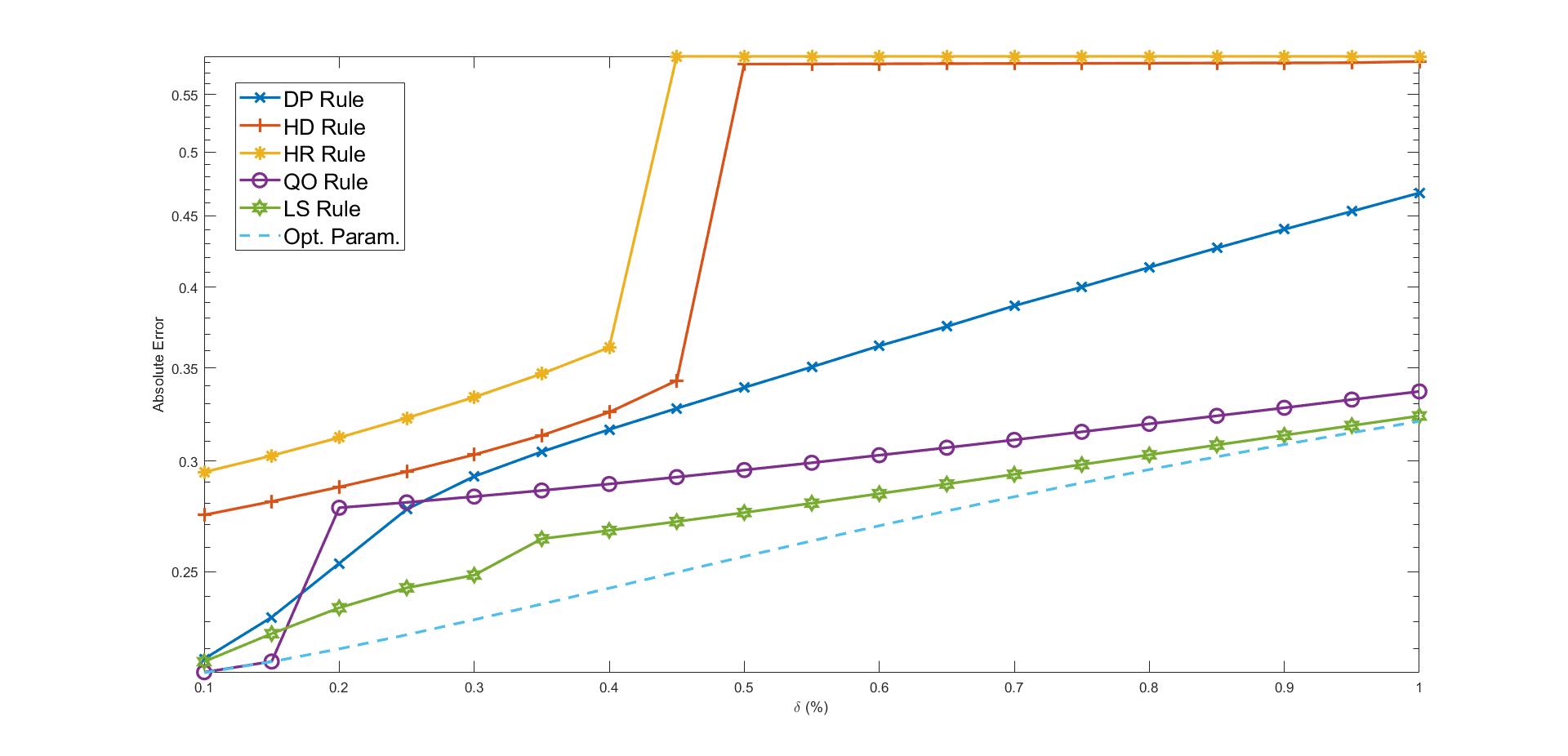}
    \caption{Numerical results for the diffusion-coefficient estimation problem introduced in Section~\ref{subsect_diffusion}: Absolute error \eqref{abs_error} at the stopping indices $\ks$ determined by the discrepancy principle \eqref{discrepancy_principle}, the different heuristic parameter choice rules \eqref{heuristic_rules}, and the optimal stopping index $\kopt$ defined in \eqref{def_kopt}, each for different relative noise levels $\delta$.}
    \label{fig_Diffusion_results}
\end{figure}

\begin{figure}[ht!]
    \centering
    \includegraphics[width=\textwidth, trim = {6.5cm 1.5cm 6cm 2cm}, clip = true]{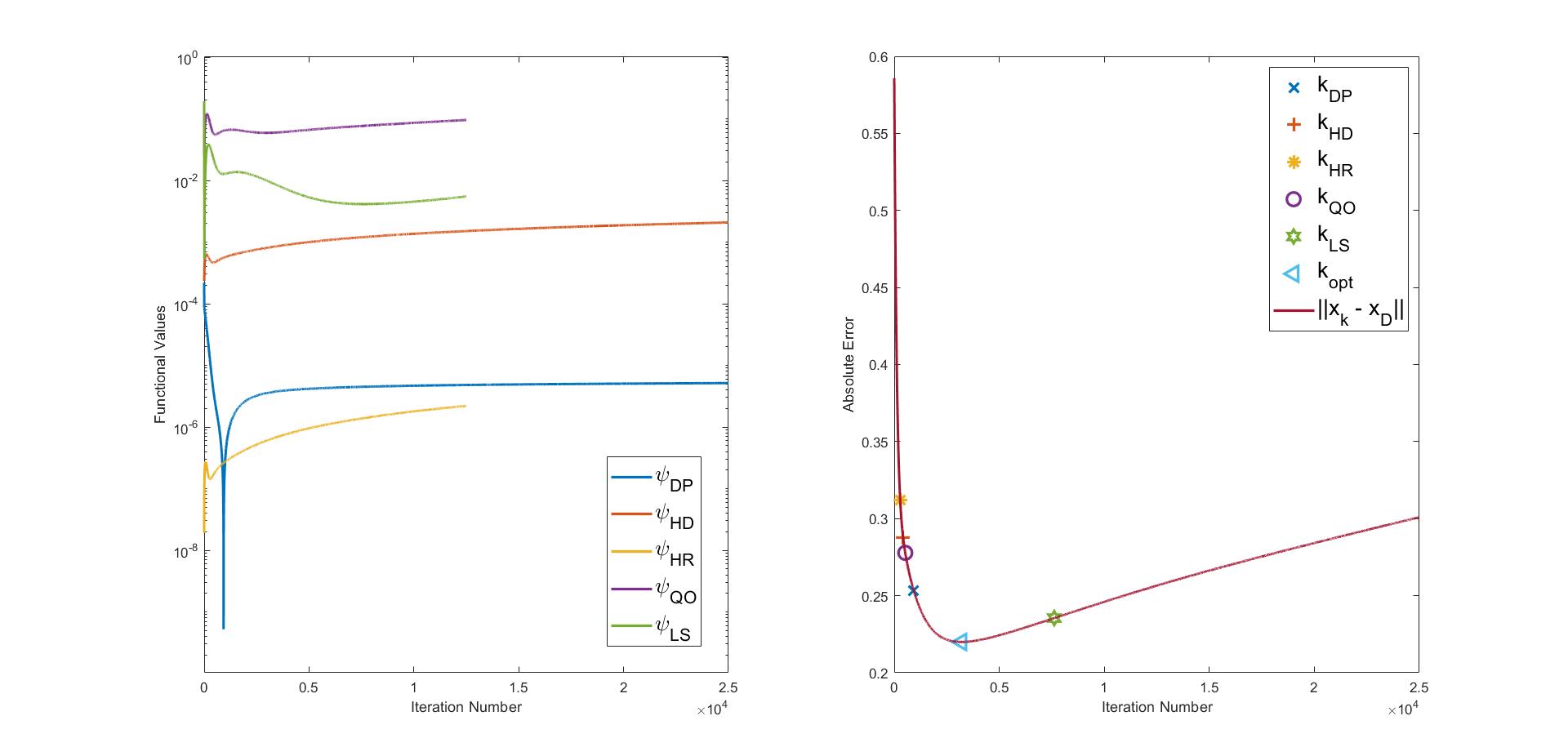}
    \caption{Numerical results for the diffusion-coefficient estimation problem introduced in Section~\ref{subsect_diffusion}: Heuristic functionals \eqref{heuristic_rules} and discrepancy functional \eqref{def_hrDP} for $\delta = 0.2\%$ relative noise (left). Corresponding evolution of the absolute error $\norm{\xkd - \xd}$ with marked points indicating the stopping indices chosen by the different rules (right).
    }
    \label{fig_Diffusion_functionals}
\end{figure}

Next, we consider the diffusion-coefficient estimation problem introduced in Section~\ref{subsect_diffusion}. For discretizing the problem we use a standard projection approach (see e.g.\ \cite{Engl_Hanke_Neubauer_1996}) onto a finite-dimensional subspace of $H^1[0,1]$ spanned by piecewise linear FEM hat functions defined on a uniform subdivision of $[0,1]$ into $50$ subintervals. For the exact solution we choose $\xd(s) = 2 + s(1-s)$ and compute the corresponding data $y$ by applying the operator $F$ defined in \eqref{def_diffusion}, using a finer grid in order to avoid an inverse crime. For the initial guess we use $x_0(s) = 2.1$, and in the discrepancy principle \eqref{discrepancy_principle} we choose $\tau = 1.1$. As before, Figure~\ref{fig_Diffusion_results} depicts the absolute errors \eqref{abs_error} corresponding to the different parameter choice rules, now for noise levels $\delta$ between $0.1\%$ and $1\%$. A typical plot of the heuristic functionals $\hr$ and the evolution of the absolute error can be seen in Figure~\ref{fig_Diffusion_functionals}

We observe again a slight superiority of the QO and LS rules over the HD and HR methods and even over the discrepancy principle, although the rate of these methods (the slope of the plots in Figure~\ref{fig_Diffusion_results})  are comparable. The HD and HR method indicate erratic behaviour for large $\delta$, which is explained by a lack of a clear minimum and the resulting sub-optimal choice of the spurious minimum. While the discrepancy principle follows quite nicely a theoretically predicted rate, the QO rule (and less the LS rule) has a jump in the error curve, which is explained by the occurrence of two local minima in the graph of Figure~\ref{fig_Diffusion_functionals}; at a certain $\delta$ the global minimum switches from one to the other. As before, we can conclude sucessful results for the discrepancy principle as well as for the heuristic rules.

\subsection{Acousto-Electrical Tomography}

Next, we consider the AET problem introduced in Section~\ref{subsect_AET}. For a detailed description of the problem setup, discretization, and implementation, we refer to \cite{Hubmer_Knudsen_Li_Sherina_2018}. In short, the unknown inclusion $\sigma^\dagger$ consists of three uniform disconnected inclusions (two circular, one crescent shaped) with values of $1.3$, $1.7$, and $2$, respectively, in an otherwise constant background of value $1$ over the circular domain $\Omega := \Kl{(r,\theta) \in [0,1) \times [0,2\pi]} \subset \R^2$. Furthermore, we use the boundary flux functions 
    \begin{equation*}
        g_j(r,\theta) := 
        \begin{cases}
            \sin(2j\pi\theta/\alpha) \,, 
            & (r,\theta) \in \Gamma(\alpha) \,,
            \\
            0 \,, & \text{else} \,,
        \end{cases}
        \qquad
        \text{for} \, j = 1,2,3 \,,
    \end{equation*}
where $\Gamma(\alpha) := \Kl{(r,\theta) \in \Kl{1} \times [0,\alpha]} \subset \partial \Omega$ for $\alpha \in [0,2\pi]$. Hence, if $\alpha = 2\pi$ then $g_j \neq 0$ almost everywhere on $\partial \Omega$. In the following, this case will be called \emph{$100\%$ boundary data}, while the case $\alpha = 3\pi/2$ is analogously called \emph{$75\%$ boundary data}. Note that in the $75\%$ boundary case, the inverse problem shows only mild instability, while in the $100\%$ case the problem behaves essentially like a well-posed problem. This can be quantified via the condition number of the discretized Fr\'echet derivative of the underlying nonlinear operator, which is equal to $385$ and $12$ in the $75\%$ and $100\%$ boundary data cases, respectively. The power density data $E_j(\sigma^\dagger)$ is created via solving the PDE \eqref{eq_AET_PDE} for $\sigma = \sigma^\dagger$, and for the initial guess we use $\sigma_0(r,\theta) = 1.5$. For completeness, note that on the definition space of the underlying nonlinear operator we use the same weighted inner product as described in \cite{Hubmer_Knudsen_Li_Sherina_2018}.

\begin{figure}[ht!]
    \centering
    \includegraphics[width=\textwidth, trim = {6.5cm 1.5cm 6cm 2cm}, clip = true]{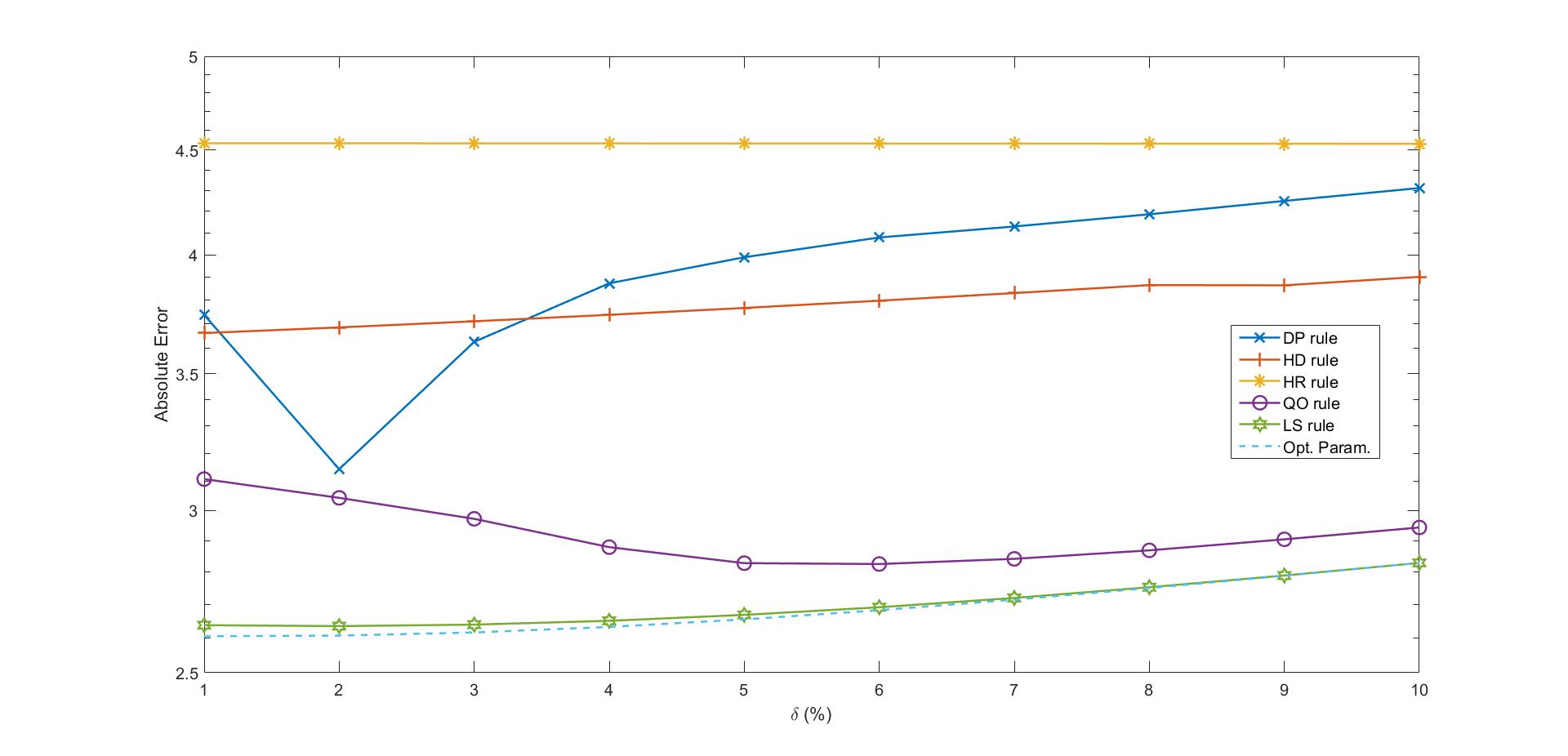}
    \caption{Numerical results for the AET problem introduced in Section~\ref{subsect_AET} with $75\%$ boundary data: Absolute error \eqref{abs_error} at the stopping indices $\ks$ determined by the discrepancy principle \eqref{discrepancy_principle}, the different heuristic parameter choice rules \eqref{heuristic_rules}, and the optimal stopping index $\kopt$ defined in \eqref{def_kopt}, each for different relative noise levels $\delta$.}
    \label{fig_AET_results_75}
\end{figure}

\begin{figure}[ht!]
    \centering
    \includegraphics[width=\textwidth, trim = {6.5cm 1.5cm 6cm 2cm}, clip = true]{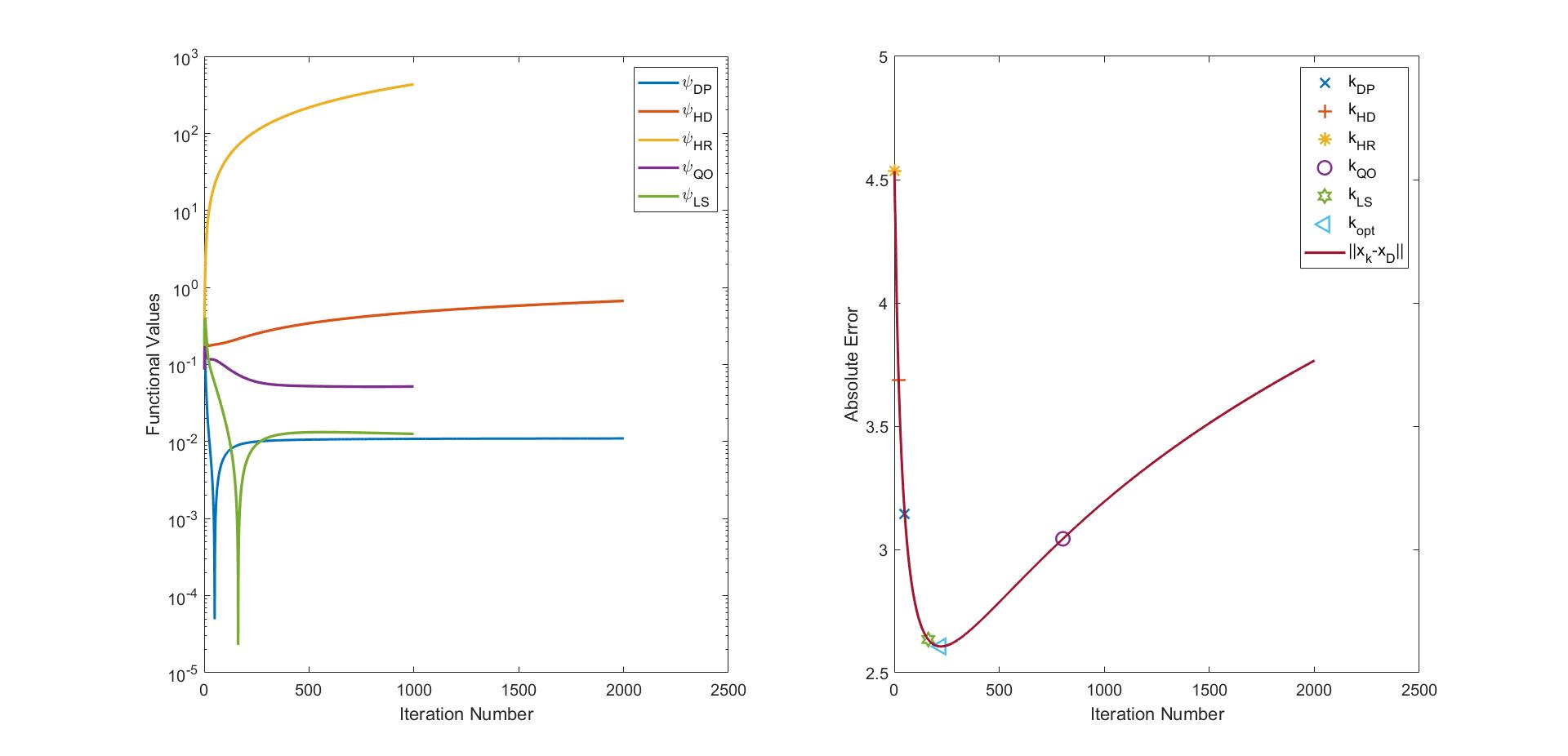}
    \caption{Numerical results for the AET problem introduced in Section~\ref{subsect_AET} with $75\%$ boundary data: Heuristic functionals \eqref{heuristic_rules} and discrepancy functional \eqref{def_hrDP} for $\delta = 2\%$ relative noise (left). Corresponding evolution of the absolute error $\norm{\xkd - \xd}$ with marked points indicating the stopping indices chosen by the different rules (right).}
    \label{fig_AET_functionals_75}
\end{figure}

First, we present results for the $75\%$ boundary data case, for which the resulting absolute errors \eqref{abs_error} corresponding to the different parameter choice rules can be found in Figures~\ref{fig_AET_results_75}. As previously, characteristic plots of the heuristic functionals $\hr$ and the evolution of the absolute errors can be found in Figure~\ref{fig_AET_functionals_75}. First of all, consider the results of discrepancy principle \eqref{discrepancy_principle}, here used with $\tau = 2$. While its corresponding stopping index is not too far off the optimal value, the steep shape of the absolute error curve nevertheless results in an overall large error. While this suggests to use a smaller $\tau$, note that already in this case for $\delta = 1\%$ the discrepancy principle is not attainable within a reasonable number of iterations. Next, note that the HR rule fails, since its corresponding functional $\hr_\text{HR}$ is monotonously increasing. From the remaining heuristic parameter choice rules, the LS rule gives the best results overall, determining a stopping index close to the optimal one. In contrast, the HD and QO role stop the iteration relatively early and late, respectively, and thus lead to suboptimal absolute errors.

\begin{figure}[ht!]
    \centering
    \includegraphics[width=\textwidth, trim = {6.5cm 1.5cm 6cm 2cm}, clip = true]{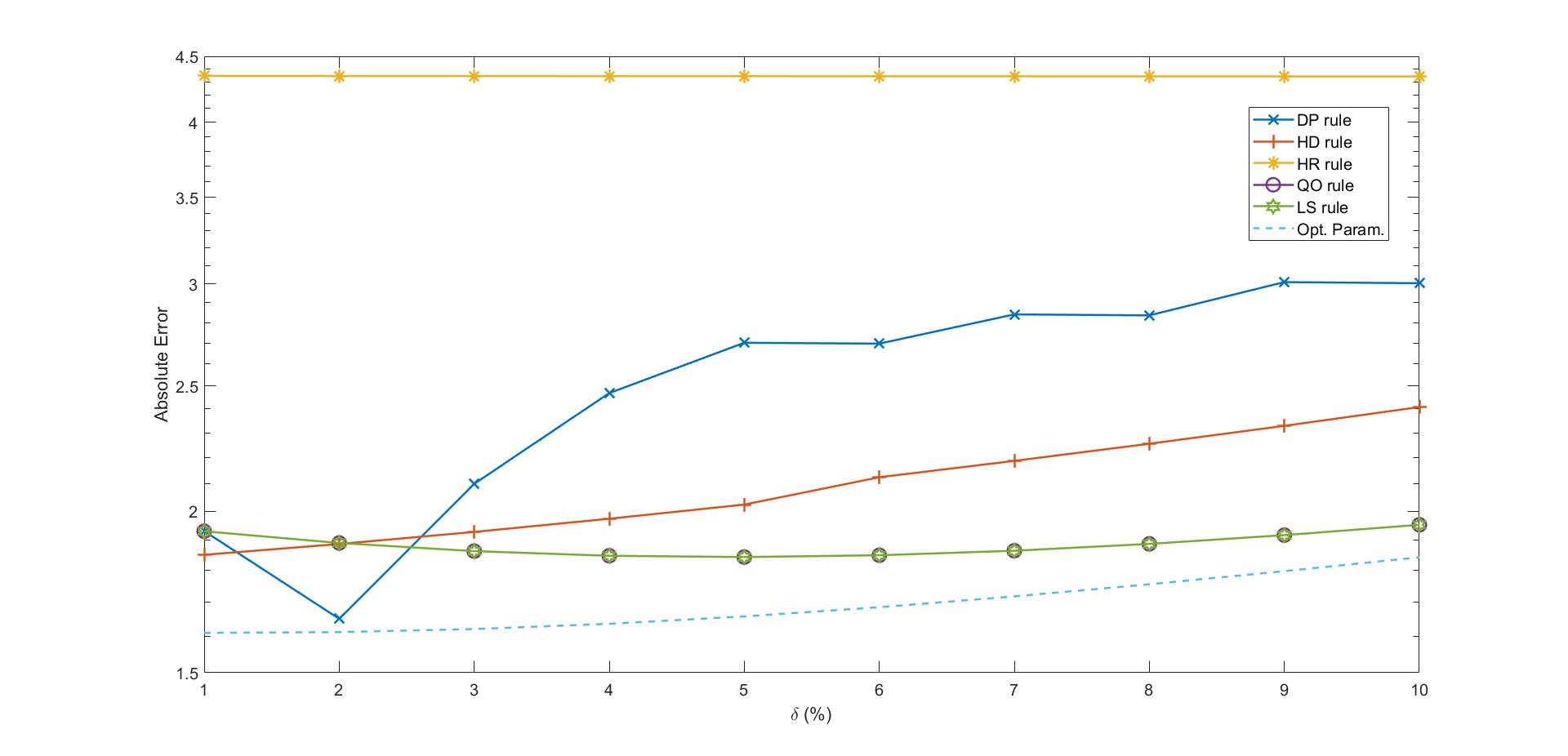}
    \caption{Numerical results for the AET problem introduced in Section~\ref{subsect_AET} with $100\%$ boundary data: Absolute error \eqref{abs_error} at the stopping indices $\ks$ determined by the discrepancy principle \eqref{discrepancy_principle}, the different heuristic parameter choice rules \eqref{heuristic_rules}, and the optimal stopping index $\kopt$ defined in \eqref{def_kopt}, each for different relative noise levels $\delta$.}
    \label{fig_AET_results_100}
\end{figure}

\begin{figure}[ht!]
    \centering
    \includegraphics[width=\textwidth, trim = {6.5cm 1.5cm 6cm 2cm}, clip = true]{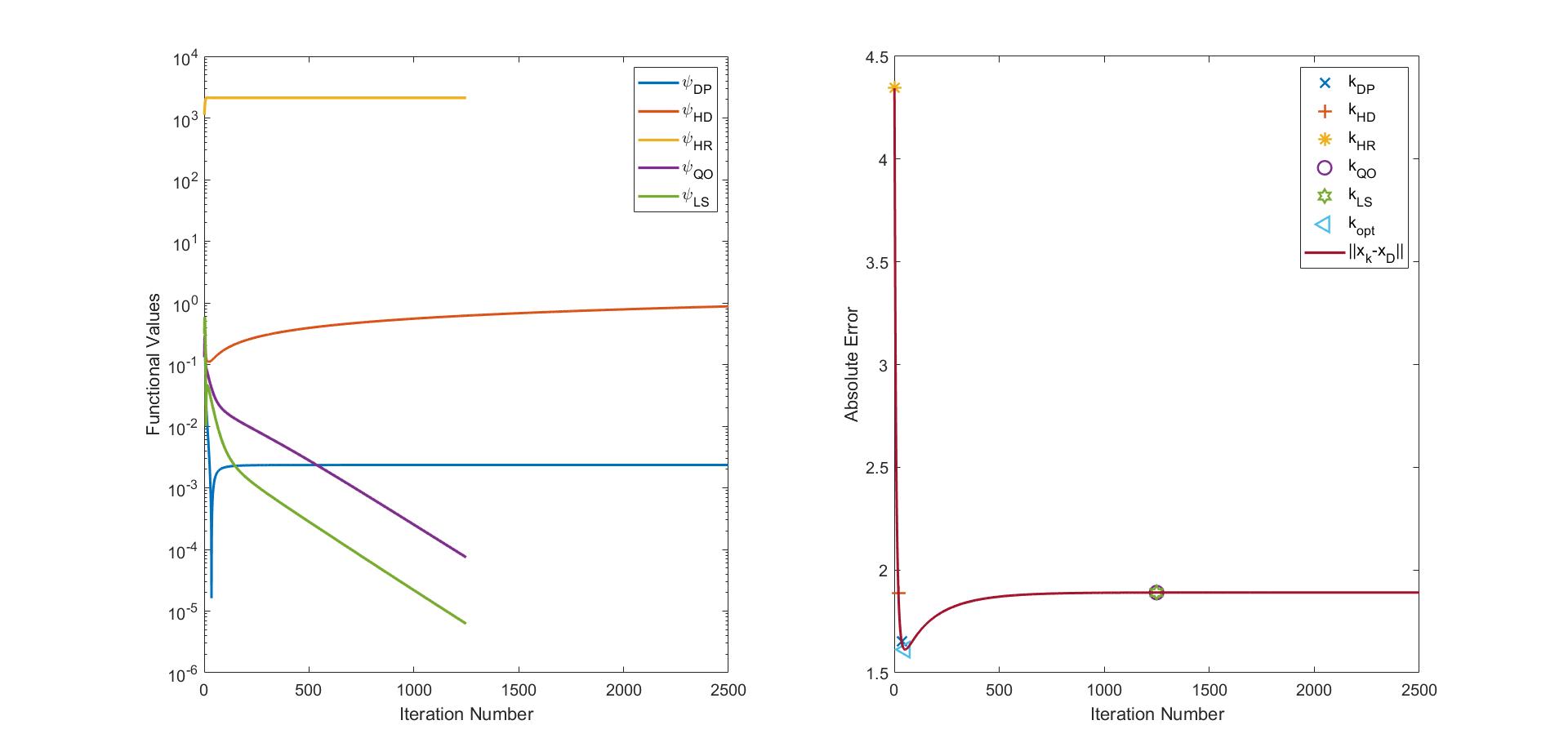}
    \caption{Numerical results for the AET problem introduced in Section~\ref{subsect_AET} with $100\%$ boundary data: Heuristic functionals \eqref{heuristic_rules} and discrepancy functional \eqref{def_hrDP} for $\delta = 2\%$ relative noise (left). Corresponding evolution of the absolute error $\norm{\xkd - \xd}$ with marked points indicating the stopping indices chosen by the different rules (right).}
    \label{fig_AET_functionals_100}
\end{figure}

Next, we consider the results for the $100\%$ boundary data case, for which the resulting absolute errors \eqref{abs_error} corresponding to the different parameter choice rules can be found in Figures~\ref{fig_AET_results_100}. Characteristic plots of the heuristic functionals $\hr$ and the evolution of the absolute errors can now be found in Figure~\ref{fig_AET_functionals_100}. Following our previous findings, we now consider the discrepancy principle \eqref{discrepancy_principle} with the choice $\tau = 1.1$. However, apart from the case of $\delta = 2\%$ the resulting stopping index is still far away from the optimal one, with the case of $\delta = 1\%$ being non-attainable as before. Next, note that both the HR, QO, and LS rule fail, with their corresponding functionals being either monotonously increasing or decreasing. We conjecture that this failure 
might be  related to the fact that the Muckenhoupt condition
is not satisfied due to a nearly well-posed situation; 
 cf.~Remark~\ref{rem_heuremark}. 
As noted above, the $100\%$ boundary case behaves essentially like a well-posed problem, which is reflected in the evolution of the absolute error depicted in Figure~\ref{fig_AET_functionals_100} (right).
Consequently a suboptimal stopping index has little effect, and 
the resulting errors when using the QO and the LS rule are comparable to those of the HD rule, which in this setting is the only heuristic rule producing a well-defined stopping index.

\subsection{SPECT}

\begin{figure}[ht!]
    \centering
    \includegraphics[width=\textwidth, trim = {6.5cm 1.5cm 6cm 2cm}, clip = true]{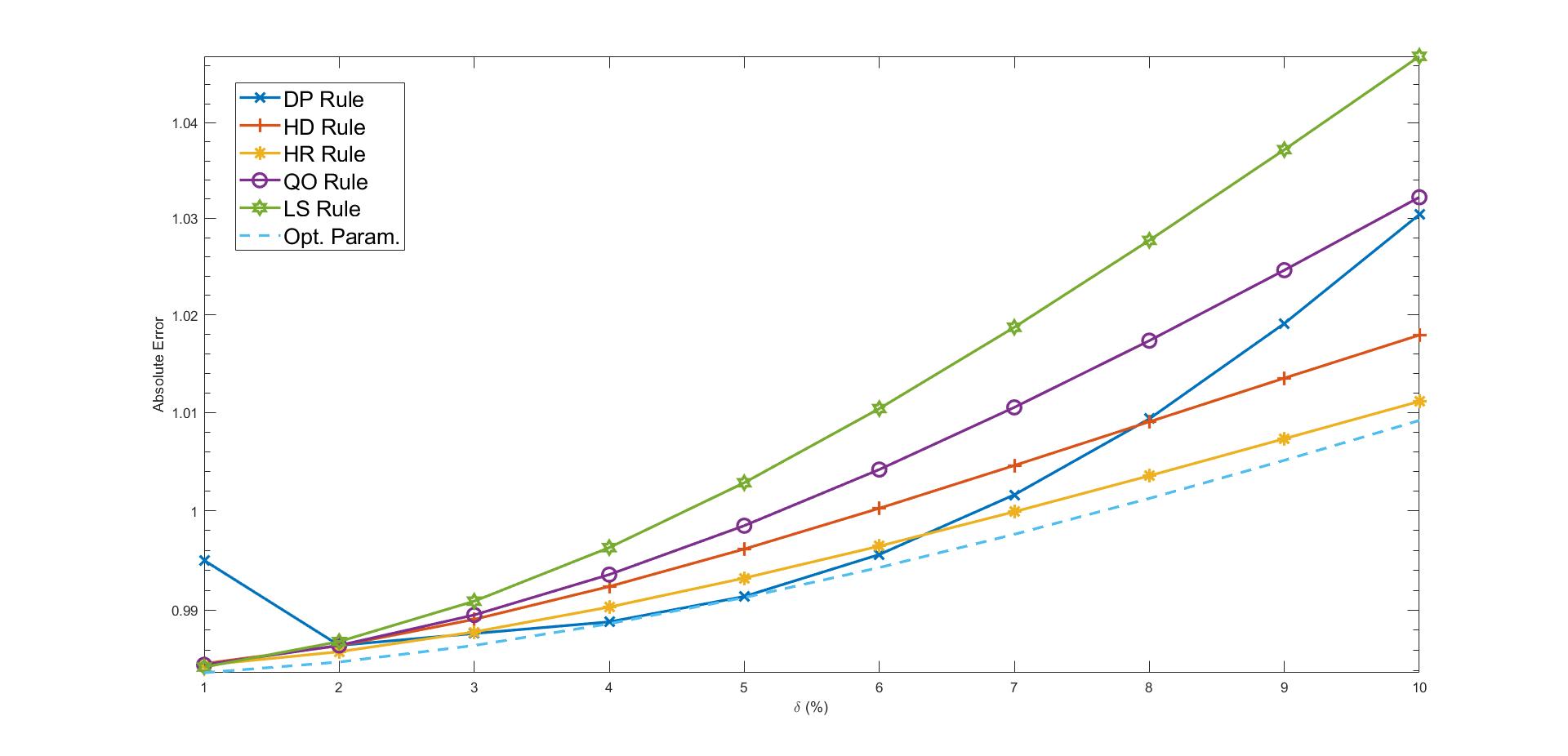}
    \caption{Numerical results for the nonlinear SPECT problem introduced in Section~\ref{subsect_SPECT}: Absolute error \eqref{abs_error} at the stopping indices $\ks$ determined by the discrepancy principle \eqref{discrepancy_principle}, the different heuristic parameter choice rules \eqref{heuristic_rules}, and the optimal stopping index $\kopt$ defined in \eqref{def_kopt}, each for different relative noise levels $\delta$.}
    \label{fig_SPECT_results}
\end{figure}

\begin{figure}[ht!]
    \centering
    \includegraphics[width=\textwidth, trim = {6.5cm 1.5cm 6cm 2cm}, clip = true]{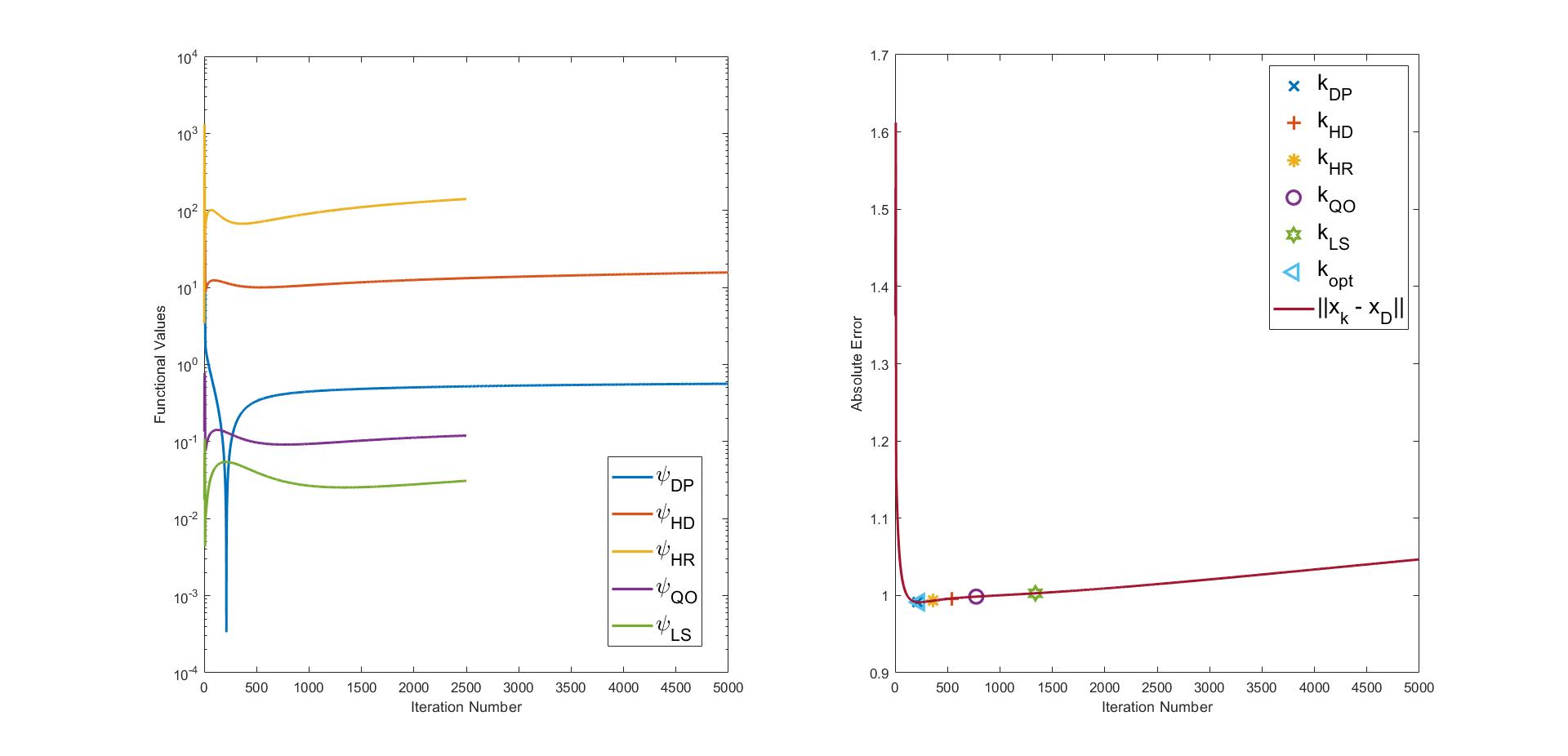}
    \caption{Numerical results for the nonlinear SPECT problem introduced in Section~\ref{subsect_SPECT}: Heuristic functionals \eqref{heuristic_rules} and discrepancy functional \eqref{def_hrDP} for $\delta = 5\%$ relative noise (left). Corresponding evolution of the absolute error $\norm{\xkd - \xd}$ with marked points indicating the stopping indices chosen by the different rules (right).}
    \label{fig_SPECT_functionals}
\end{figure}

For the fifth test, we consider the nonlinear SPECT problem introduced in Section~\ref{subsect_SPECT}. For discretizing the problem we utilize the same approach used e.g.\ in \cite{Hubmer_Ramlau_2017,Ramlau_2003,Ramlau_Teschke_2006}, using $79$ uniformly spaced angles $\omega$ in the attenuated Radon transform \eqref{def_SPECT}. For the exact solution $(f^\dagger,\mu^\dagger)$, we choose the MCAT phantom \cite{Terry_Tsui_Perry_Hendricks_Gullberg_1990}, see also \cite{Hubmer_Ramlau_2017,Ramlau_2003,Ramlau_Teschke_2006}, and for the initial guess we use $(f_0,\mu_0) = (0,0)$. For the discrepancy principle \eqref{discrepancy_principle} we choose the rather large value $\tau = 10$, which was however found to lead to better results than the standard choices $\tau = 1.1$ or $\tau = 2$. The absolute error \eqref{abs_error} corresponding to the different parameter choice rules are depicted in Figure~\ref{fig_SPECT_results}, in this case for the practically realsitc case of noise levels $\delta$ between $1\%$ and $10\%$. Again, typical plots of the heuristic functionals $\hr$ as well as the evolution of the absolute error over the iteration is depicted in Figure~\ref{fig_SPECT_functionals}.

The plots in these figures show that all parameter choice rules work well for this problem, that the error rate follow the optimal one, and that the heuristic plots in Figure~\ref{fig_SPECT_functionals}
have a clear minimum. Furthermore, the HD and HR rules are slightly superior to the other rules in terms of the resulting absolute error. Note that as for the nonlinear Hammerstein problem, the error graph for the discrepancy principle again exhibits a slight parabola shape, which could again be related either to a too small choice of $\tau$, or to a coarse discretization of the problem.

\subsection{Auto-Convolution}

\begin{figure}[ht!]
    \centering
    \includegraphics[width=\textwidth, trim = {6.5cm 1.5cm 6cm 2cm}, clip = true]{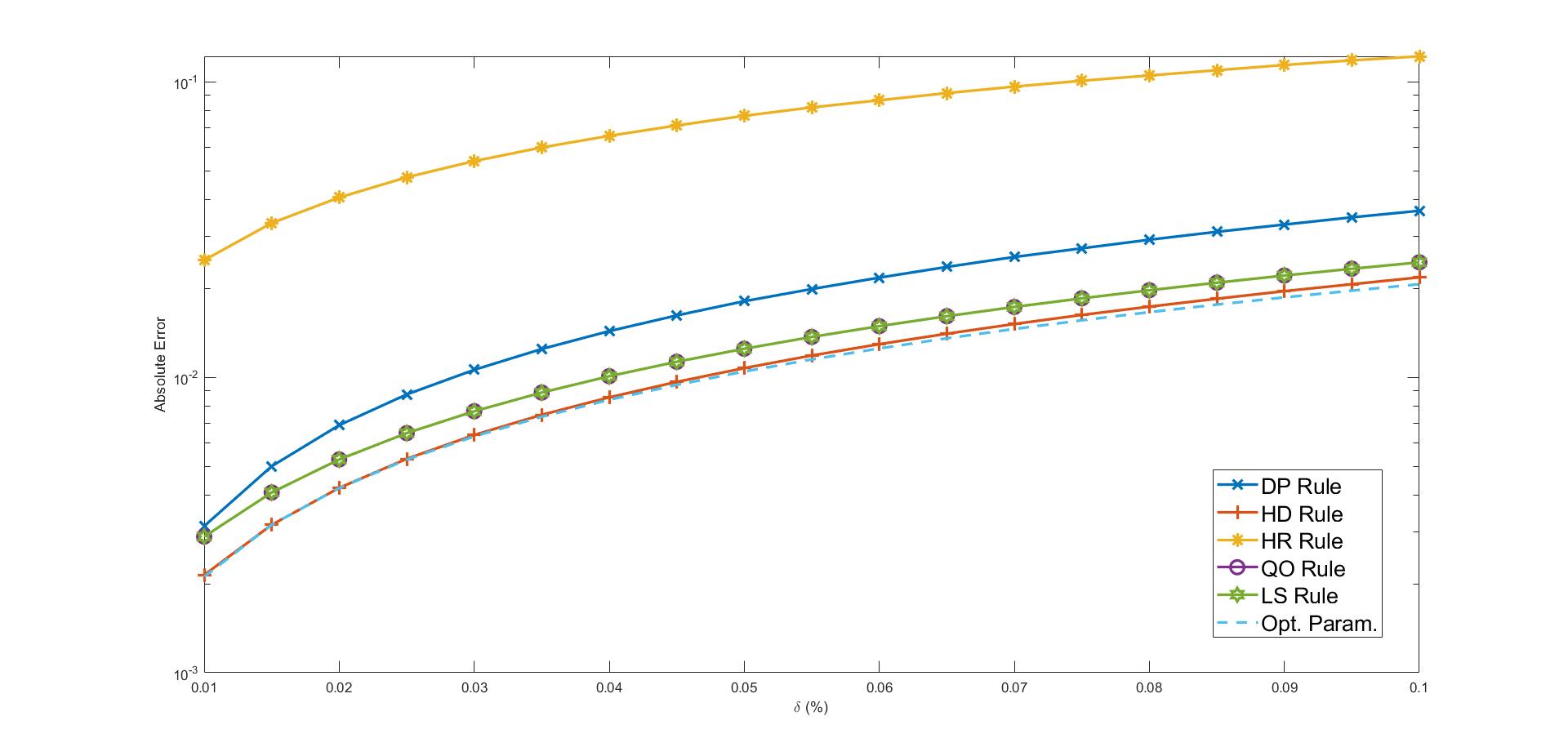}
    \caption{Numerical results for the auto-convolution problem introduced in Section~\ref{subsect_conv}: Absolute error \eqref{abs_error} at the stopping indices $\ks$ determined by the discrepancy principle \eqref{discrepancy_principle}, the different heuristic parameter choice rules \eqref{heuristic_rules}, and the optimal stopping index $\kopt$ defined in \eqref{def_kopt}, each for different relative noise levels $\delta$.}
    \label{fig_Autoconvolution_results}
\end{figure}

\begin{figure}[ht!]
    \centering
    \includegraphics[width=\textwidth, trim = {6.5cm 1.5cm 6cm 2cm}, clip = true]{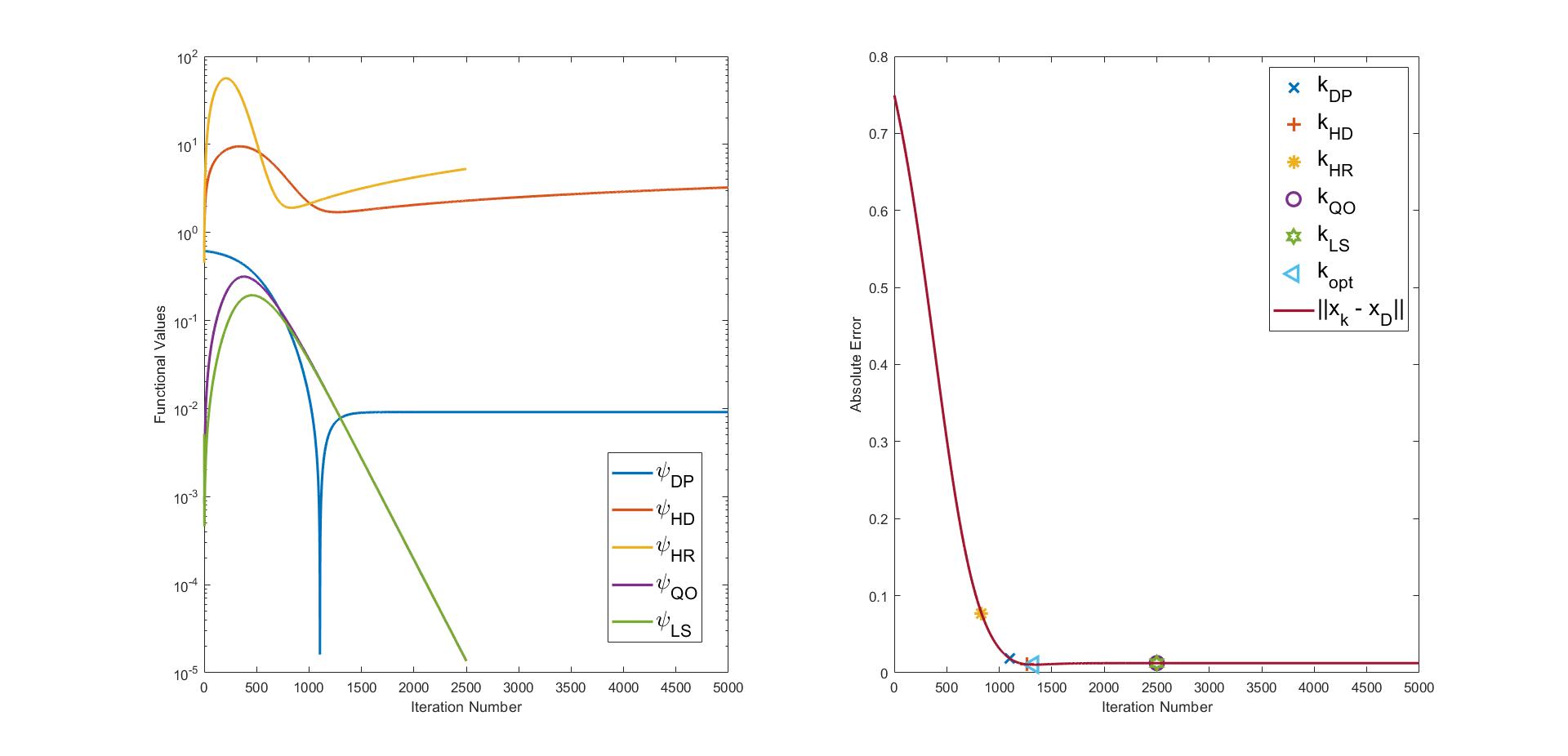}
    \caption{Numerical results for the auto-convolution problem introduced in Section~\ref{subsect_conv}: Heuristic functionals \eqref{heuristic_rules} and discrepancy functional \eqref{def_hrDP} for $\delta = 0.05\%$ relative noise (left). Corresponding evolution of the absolute error $\norm{\xkd - \xd}$ with marked points indicating the stopping indices chosen by the different rules (right).}
    \label{fig_Autoconvolution_functionals}
\end{figure}

For the final test, we consider the auto-convolution problem introduced in Section~\ref{subsect_conv}. For the discretization of this problem, which is based on standard FEM hat functions on a uniform subdivision of the interval $[0,1]$ into $60$ subintervals, we refer to \cite{Hubmer_Ramlau_2018}. For the exact solution we choose $\xd(s) = 10 + \sqrt{2} \sin(2 \pi s)$, from which we compute the corresponding data $y$ by applying the operator $F$ as defined in \eqref{def_conv}. For the initial guess we use $x_0(s) = 10 + \tfrac{1}{4}\sqrt{2} \sin(2 \pi s)$, and in the discrepancy principle we choose $\tau = 1.1$. Since the initial guess is rather close to the exact solution, we now consider noise levels $\delta$ between $0.01 \%$ and $0.1\%$. The corresponding absolute errors \eqref{abs_error} for the different parameter choice rules are depicted in Figure~\ref{fig_Autoconvolution_results}, while a typical plot of the heuristic functionals $\hr$ and the evolution of the absolute error can be seen in Figure~\ref{fig_Autoconvolution_functionals}.

We observe that for this problem the HD rules give the best results overall, yielding stopping indices close to the optimal $\kopt$ for all considered noise levels. Furthermore, from Figure~\ref{fig_Autoconvolution_functionals} we can see that also the functional $\hr_\text{HR}$ exhibits a clearly distinguishable minimum. However, since the absolute error drops steeply between the corresponding stopping index and the optimal one, the resulting absolute error when using the HR rule is significantly higher. In contrast, both the QO rule and the LS rule tend towards $-\infty$ as the iteration number increases. This coincides with the fact that the absolute error in the iteration stays more or less constant after having reached the minimum value at $\kopt$. Consequently, the LS and the QR rule both stop with $\ks = \kmax$, which in this case by chance leads to generally very good absolute errors. However, the shape of the QO and LS functionals is in contrast to theory, as the graph is expected to diverge for $k\to\infty$. This observation is a hint that the Muckenhoupt condition is not satisfied here for QO and LS, while it might be for HD and HR; cf.~Remark~\ref{rem_heuremark}. This may also be connected to the fact that the residual functional $x \mapsto \norm{F(x)-\yd}^2$ is locally convex around our chosen solution $\xd$, see e.g.\ \cite[Proposition~5.2]{Hubmer_Ramlau_2018}, and thus the auto-convolution problem may behave nearly like a well-posed problem; cf.~Section~\ref{subsect_challenges}.

\section{Summary and Conclusion}\label{sect_conclusion}

In the previous section, we presented and discussed the results of the different parameter choice rules applied to the test problems introduced in Section~\ref{sect_test_problems}. As a summary, we now present our findings in Table~\ref{tabtab}, where we classify the results, in a rather informal style, together with comments and suspected background. 

\begin{table}[ht!]
\caption{Summary of performance for various examples and stopping rules.}\label{tabtab}
\begin{center}
\begin{tabular}{|m{3.3cm}|c||m{0.9cm}|m{1.5cm}|m{1.5cm}|m{1.5cm}|m{1.5cm}|}
    \multicolumn{1}{c|}{Example} &Indicator &$\delta$-rule  &\multicolumn{4}{c}{Heuristic Rules} \\ \hline &  &\rule{0mm}{2.3ex}DP & HD & HR & QO & LS  \\ \hline \hline 
    \multirow{2}{*}{\makecell{\\[-7mm] Hammerstein\\ equation \\ Illposedness: Mild \\ TCC: Yes }} & $k_*\sim k_{opt}$ & \Go & \multicolumn{2}{c|}{\makecell{$\delta$ small: \Go \\ $\delta$ large: \Ba }} & \multicolumn{2}{c|}{\Av} \\ \cline{2-7} & Error: & \Go & \multicolumn{2}{l|}{\makecell{$\delta$ small: \Go \\ $\delta$ large: \Ba }} & \multicolumn{2}{c|}{\Av} \\ \hline \hline
    \multirow{2}{*}{\makecell{\\[-7mm] Diffusion\\ Estimation \\ Illposedness: Mild \\ TCC: Yes }} & $k_*\sim k_{opt}$ & \Go & \multicolumn{2}{c|}{\makecell{$\delta$ small: \Go \\ $\delta$ large: \Ba }} & \multicolumn{2}{c|}{\Go} \\ \cline{2-7} & Error: & \Go & \multicolumn{2}{l|}{\makecell{$\delta$ small: \Go \\ $\delta$ large: \Ba }} & \multicolumn{2}{c|}{\Go} \\ \hline \hline 
     \multirow{2}{*}{\makecell{\\[-3mm] Acousto-Electric \\ Tomography $75\%$ \\ Illposedness: Mild \\ TCC: Unknown }} & $k_*\sim k_{opt}$ & \makecell{ \Go \\ } & \multicolumn{2}{c|}{\Ba} &  \Av & \Ex  \\  \cline{2-7} & Error: & \makecell{ \Go \\ } & \Go & \Go & \Go & \Ex \\ \hline \hline 
     \multirow{2}{*}{\makecell{\\[-3mm] Acousto-Electric\\ Tomography $100\%$ \\ Near-Wellposed \\ TCC: Unknown }} & $k_*\sim k_{opt}$ & \makecell{ \Go \\ } & \Go & \multicolumn{3}{c|}{\Ba} \\ \cline{2-7} & Error: &  \makecell{ \Go \\ } & \Go & \Av & \Go & \Ex \\ \hline \hline
    \multirow{2}{*}{\makecell{\\[-6mm] \\ SPECT \\ Illposed: Mild \\ TCC: Unknown }} & $k_*\sim k_{opt}$ & \makecell{ \Go \\ } & \multicolumn{4}{c|}{\Go} \\ \cline{2-7} & Error: & \makecell{ \Go \\ } & \Go & \Ex & \Go & \Go \\ \hline \hline
    \multirow{2}{*}{\makecell{\\[-6mm]  Auto-\\ Convolution \\ Illposed: Mild \\ TCC: Unknown }} & $k_*\sim k_{opt}$ & \makecell{ \Go \\ } & \multicolumn{2}{c|}{\Go} & \multicolumn{2}{c|}{\Ba} \\ \cline{2-7} & Error: & \makecell{ \Go \\ } & \Ex & \Av & \Av & \Ex \\ \hline 
\end{tabular} 
\end{center} 
\end{table}

In the first column of Table~\ref{tabtab}, we indicate the (suspected) type of ill-posedness of the various test problems and whether the tangential cone condition (TCC) is known to hold. Furthermore, in the second column the performance of the parameter choice rules is classified according to two indicators: The row with $k_* \sim \kopt$ indicates whether the stopping index selected by the respective parameter choice rule is close to the optimal index, and whether the heuristic functionals behave as desired, i.e., have a clear minimum. The row with ``Error'' indicates how close the resulting absolute error is to the optimal error. These performance indicators are stated in a colloquial manner, classified as ``Excellent'', ``Good'', ``Average'', or ``Bad''.  

Some conclusions can be drawn from these results: First of all, the discrepancy principle works well in all cases, even when the tangential cone condition is not known to hold; as discussed above. However, it requires a proper choice of the parameter $\tau$, and the strange parabola shape in Figure~\ref{fig_Hammerstein_results} is, e.g., attributed to its having been selected too small. Next, we observe that the heuristic rules work well for many cases, but not always. Furthermore, our numerical results show that it is not possible to determine an overall ``best'' heuristic parameter choice rule for all test problems. This suggests that in practice one should always first conduct a series of simulations using multiple different parameter choice rules for any given problem, instead of blindly selecting any single rule among them. From the computational point of view, the authors personally prefer the HD rule, which does not require to compute and store both iterates $x_k$ and $x_{2k}$ during the iteration. However, also the HD rule fails at times, and the LS rule in particular has been found to be a useful alternative which often showed an excellent performance with respect to the resulting absolute error. One issue with the QO and LS rules is that their corresponding functionals $\hr$ do not exhibit a clear minimium in the case of well-posed or nearly well-posed problems (e.g., AET with $100\%$ boundary data and Auto-Convolution). However, this can be attributed to the fact that the Muckenhoupt condition is stronger for these rules and hence might not be satisfied; cf.~Remark~\ref{rem_heuremark}. On the other hand, even if the QO and LS rule fail in finding a good approximation of $k_{opt}$, the resulting error is then not too dramatic, since the almost well-posedness then only leads to moderate errors. 

Finally, our numerical studies also have implications for the further analytic study and development of (novel) heuristic stopping rules for nonlinear Landweber iteration. In particular, our findings indicate that apart from the relation between the smoothness of the solution and the noise (i.e., noise conditions), also the type of ill-posedness and especially the type of nonlinearity of the problem have to be taken into account. The authors believe that it will in general be impossible to design and analyse a single heuristic parameter choice rules which performs well for all different nonlinear inverse problems. Rather, one may have to consider different heuristic rules depending on the type of nonlinearity and ill-posedness of each specific problem. This then poses an interesting challenge for future research.

\section{Support}
The authors were funded by the Austrian Science Fund (FWF): F6805-N36 (SH), F6807-N36 (ES), and P30157-N31 (SK and KR).

\bibliographystyle{plain}
{\footnotesize
\bibliography{mybib}
}

\end{document}